\documentclass{article}

\usepackage{arxiv}
\usepackage{amsmath,amssymb}
\usepackage[english]{babel}
\usepackage{subfigure}
\usepackage{caption}
\usepackage{array}
\usepackage{float}
\usepackage{color,soul}
\usepackage[utf8]{inputenc} % allow utf-8 input
\usepackage[T1]{fontenc}    % use 8-bit T1 fonts
\usepackage{url}            % simple URL typesetting
\usepackage{booktabs}       % professional-quality tables
\usepackage{amsfonts}       % blackboard math symbols
\usepackage{nicefrac}       % compact symbols for 1/2, etc.
\usepackage{microtype}      % microtypography
\usepackage{lipsum}		% Can be removed after putting your text content
\usepackage{graphicx}
\usepackage{doi}

\numberwithin{equation}{section}

\title{Designing Social Distancing Policies for the COVID-19 Pandemic: A probabilistic model predictive control approach}

\author{\hspace{1mm} Antonios Armaou \thanks{corresponding author} \\
	Dept. of Chemical Engineering\\
	The Pennsylvania State University, USA
	\texttt{armaou@psu.edu} \\
	\And
	\hspace{1mm}Bryce Katch \\
		Dept. of Chemical Engineering\\
	The Pennsylvania State University, USA
	\texttt{buk71@psu.edu} \\
		\And
	\hspace{1mm}Lucia Russo\\
	Institute of Science and Technology for Energy and Sustainable Mobility\\
	Consiglio Nazionale delle Ricerche, Italy
	\texttt{lucia.russo@stems.cnr.it} \\
			\And
	\hspace{1mm}Constantinos Siettos \thanks{corresponding author}\\
Dipartimento di Matematica e Applicazioni ``Renato Caccioppoli"\\
Universit\`a degli Studi di Napoli Federico II, Naples, Italy\\
	\texttt{constantinos.siettos@unina.it} \\
}

\begin{document}
\maketitle

\begin{abstract}
The effective control of the COVID-19 pandemic is one the most challenging issues of nowadays. The design of optimal control policies is perplexed from a variety of social, political, economical and epidemiological factors. Here, based on epidemiological data reported in recent studies for the Italian region of Lombardy, which experienced one of the largest and most devastating outbreaks in Europe during the first wave of the pandemic, we address a probabilistic model predictive control (PMPC) approach for the modelling and the systematic study of what if scenarios of the social distancing in a retrospective analysis for the first wave of the pandemic in Lombardy. The performance of the proposed PMPC scheme was assessed based on simulations of a compartmental model that was developed to quantify the uncertainty in the level of the asymptomatic cases in the population, and the synergistic effect of social distancing in various activities, and public awareness campaign prompting people to adopt cautious behaviors to reduce the risk of disease transmission. 
The PMPC scheme takes into account the social mixing effect, i.e. the effect of the various activities in the potential transmission of the disease. The proposed approach demonstrates the utility of a PMPC approach in addressing COVID-19 transmission and implementing public relaxation policies.
\end{abstract}

% keywords can be removed
\keywords{COVID-19 pandemic \and Model Predictive Control \and Lombardy Italy \and Social Distancing \and  Social Mixing \and Google Mobility Reports}

\section{Introduction}

The COVID-19 pandemic emerged in China in early December of 2019, two years ago sweeping the globe in a short period of time. During the first wave, by May 20, 2020, a total of more than 4,900,000 cases in 188 countries were reported  while the death toll was 320,000 (\cite{Johnsmap}). Today, the death toll worldwide has passed the 5.2 million people and the total infected cases are more than 267 million. After China, Italy became the global epicenter of the pandemic from March to mid-April, listing hundreds of deaths each day and surpassing soon the fatalities of China. By May 19, the number of confirmed cases exceeded 225,000, while the death toll passed 32,000. Of all regions, Lombardy\textemdash the nation's economic center\textemdash suffered the most; more than half of the deaths counted in Italy came from Lombardy. The first confirmed case in the region was reported on February 21, 2020, and on March 8, the government announced a country-wise lockdown which was expanded on March 21 to include even stricter isolation measures such as the prohibition of mobility and closure of all non-essential business activities. According to Google Community Mobility Reports, population mobility in various activities was reduced as much as 80\% compared to the pre-pandemic baseline activity \cite{Google}. By the end of March 2020, as the disease continued to spread rapidly, the fear of the collapse of the healthcare systems around the world led more than 100 countries to implement complete or partial lockdowns.
Since the end of April 2020, due to the severe isolation measures, the number of daily new cases and deaths declined in many European and Asia-Pacific countries and governments began to introduce a phased reopening of the activities. While in most of the countries there was a decline in the confirmed cases during the summer period of 2020, a second even more intense wave of the pandemic emerged during Autumn of 2020-Winter 2021. The number of dead people due to COVID-19 at the second wave surpassed the 900 per day in France in mid November, the 500 in Spain at the end of November and the 700 early February, the 1200 in Germany at the end of December 2020 the 1500 in UK in mid of January 2021, and the 730 in Italy early December 2020 \cite{Johnsmap}. In early November 2020, in a race against time there were also the first announcements of vaccines that have been tested successfully in the first interim analysis from phase three. However, one year later, due to the heterogeneity and the generally still low rates of vaccinations worldwide, there are concerns that the COVID-19 pandemic will become seasonal as the influenza. Thus, until a global heard immunity is achieved, the systematic scheduling of social distancing measures will be still on the table of public health policies due to the relatively high mortality ratios of the disease among the elderly and vulnerable population.\par
At the beginning of the pandemic the efforts were focused on the development of mathematical models in an attempt to assess the dynamics  for forecasting purposes, as well as to estimate epidemiological parameters that cannot be computed directly based on clinical data, ranging from simple exponential growth models (see e.g. \cite{Zhao_2020, Remuzzi_2020}) and compartmental SIR-like models (\cite{Wu_2020, Anastassopoulouetal,belgaid2020analysis}) to metapopulation models \cite{Wu_2020,Chinazzi_2020} including the modelling of the influence of travel restrictions and other control measures in reducing the spread (\cite{Chinazzi_2020}).
Today, the mathematical modelling efforts are focused more on the use of such mathematical models to design in a systematic way efficient coronavirus-slowing down-control policies accounting for both the effect of vaccinations and the impact of mobility policies on the COVID-19 transmission. Simple control approaches to COVID-19 involve the use of feedback control and on/off control to simulate the implementation of mobility policies \cite{di_lauro_covid-19_2021}. More advanced control approaches such as Model predictive control (MPC) have been also  applied constrained by specific epidemiological models. For example, Kohler et al. proposed a robust MPC using a compartmental model applied to Germany \cite{kohler_robust_2020}; similar approaches have been used to investigate regions in China, Brazil, Italy, and the United States \cite{hadi_control_2021,peni_nonlinear_2020,tsay_modeling_2020,carli_model_2020}. Finally, the effect of regional interventions seen as an interacting network between regions has been also been propose \cite{della2020network}.  \par
In this work, we address a probabilistic model predictive control (PMPC) approach constrained to a compartmental model extending beyond the conventional SEIR approach to account for hospitalization and quarantine and importantly for the uncertainty in the level of asymptomatic cases in the population. In particular, the development of our model is based on epidemiological data from the region of Lombardy that have been very recently published \cite{grosso2021decreasing},  thus extending  the model proposed by \cite{Russo2020}.  In that work, \cite{Russo2020} addressed a model that includes a compartment for the asymptomatic cases and proposed an optimization algorithm to assess the uncertainty that characterizes the actual number of infected cases in the total population, which is mainly due to the large percentage of asymptomatic or mildly symptomatic cases and the uncertainty regarding the DAY-ZERO of the outbreak. The knowledge of both these quantities is crucial to assess the stage and dynamics of the epidemic, especially during the first growth period. By doing so, the proposed model was able to predict quite well the evolution of the dynamics of the pandemic in Lombardy two months ahead of time (see \cite{Russo2020}).
Here, building up on this model, we extend it by incorporating two more compartments, namely the hospitalized and the guaranteed cases as well as the effect of social mixing, i.e. the effect of various social activities in the potential transmission of infectious diseases as has been reported for Italy and other countries in \cite{mossong_social_2008}. Based on this model, we designed a PMPC scheme that exploits the use of the mobility data provided by Google Community Reports for the Lombardy region to study (retrospectively) what if scenarios of the optimal control action policies, thus comparing them with the control actions taken. Here, we should note that such retrospective analysis is a very important task as also pinpointed in a recently published study of the Italian exit strategy from COVID-19 lockdown \cite{marziano2021retrospective} as can help a better risk assessment for analogous situations in the future. 

The structure of the paper is as follows. In section 2.1, we present the compartmental dynamical model and we justify the selection of the values of the model parameters based on publicly available epidemiological data. Based on the model, in section 2.2, we derive estimations of the basic and effective reproduction numbers. In section 2.3 we perform an analysis of the Google mobility reports for Lombardy, thus providing a model for the effect of the social distancing measures in all activities to the residential activity which is also an important factor in the transmission of infectious diseases. In section 2.4, we present the proposed PMPC scheme constrained to the impact of the social mixing and activity patterns in the potential transmission. In section 3, we report the simulation results, and in section 4 we conclude with the discussion of the findings. Finally, in the Supplementary Information section, we present the sources of the data used in the present study, competing formulations to the designed PMPC one and an analysis of the effect of PMPC parameters on the resulting policies.

\section{Materials and methods}
\subsection{The Modelling Approach}\label{s:model}
As in our probabilistic model predictive control scheme, we consider as constraint a maximum number of hospitalized persons due to COVID19, we extended the model proposed in Russo et al. \cite{Russo2020} to include hospitalized cases and confirmed infected cases with mild symptoms that are isolated at home. Thus, it includes two categories of infected cases, namely the confirmed/reported and the unreported (unknown) asymptomatic cases in the total population \cite{Russo2020}. Based on observations and studies, our modelling hypothesis is that the confirmed cases of infected are only a (small) subset of the actual number of infected cases in the total population \cite{Li2020.03.06.20031880,Anastassopoulouetal}. Regarding the confirmed cases a study conducted by the Chinese CDC which was based on a total of 72,314 cases in China, about 80.9\%  of the cases were mild and could recover at home, 13.8\% severe and 4.7\%  critical \cite{CCDC11Feb}.

On the basis of the above findings, in our modelling approach, the asymptomatic cases recover from the disease relatively soon and without medical care, while the confirmed cases include all the above types, but on average their recovery lasts longer than the asymptomatic cases, they may also be hospitalized or sadly die from the disease \cite{Russo2020}.

\begin{figure}
\includegraphics[width=0.99\textwidth]{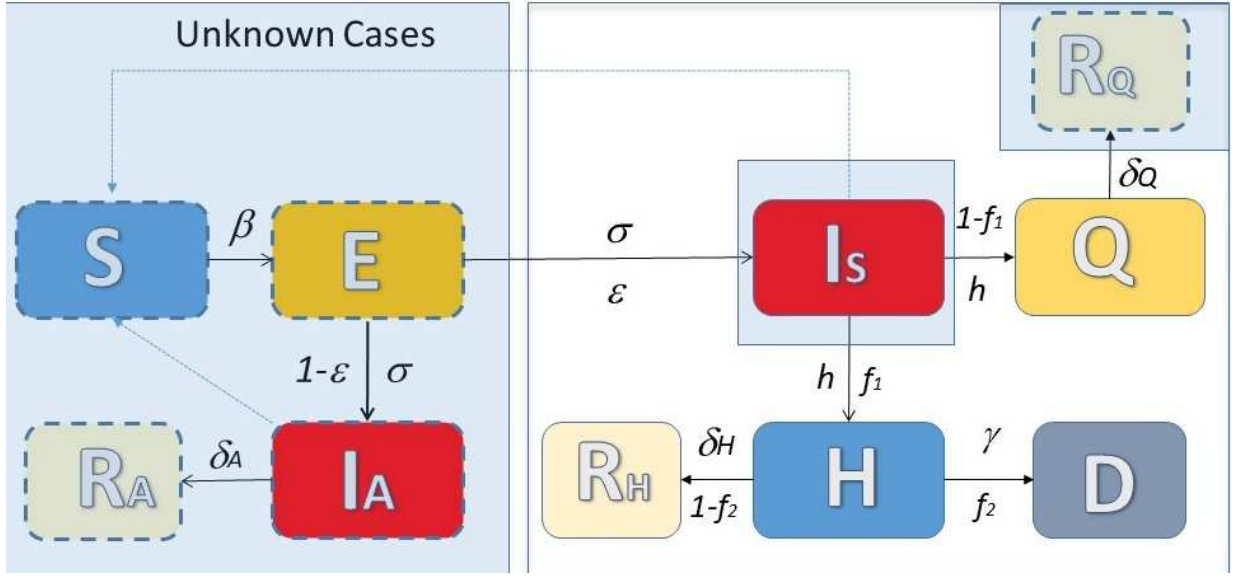}
\caption{A schematic of the proposed compartmental model with hidden compartments. Dashed lines  depict the compartments for which there are no available data. The released data report the total number of confirmed cases $I_S$, the total number of hospitalized cases $H$, the total number of recovered cases dismissed from the hospital $R_H$, the total number of confirmed infected isolated at home $Q$ and the total number of deaths $D$. 
}\label{fig0}
\end{figure}

Considering a well-mixed population of size $N$, the state of the system at time $t$ is described by the following compartments (see also Figure \ref{fig0} for a schematic): $S(t)$ representing the number of susceptible persons, $E(t)$ the number of exposed, $I_A(t)$ the number of asymptomatic infected persons in the total population that recover relatively soon without any other complications; $R(t)$ is the number of recovered asymptomatic persons. $I_S(t)$ denotes the number of infected cases who may develop more severe symptoms and a part of them is hospitalized (in our model denoted by $H(t)$), a part of them dies (denoted by $D(t)$) and a part of them recovers (denoted by $R_H(t)$). The other part of confirmed infected cases from the $I_S(t)$ compartment (in our model represented by the variable $Q(t)$) is assumed to experience mild symptoms and is quarantined at home. For this compartment there are no other available data except from its specific number. The confirmed recovered (in our model denoted by $R_H(t)$) are the recovered cases dismissed from hospitals. Furthermore, there are no reported data for the number of deaths from the infected people that are isolated at home. As said this compartment represents the part of  cases that experience mild symptoms and thus they aren't expected to die from the disease.\par
Here, we remark that a wide testing policy may also result in the identification of asymptomatic cases belonging to the compartment $I_A$ that would then be assigned to compartment $I_S$. However, as a generally reported rule in Italy, tests were conducted only for those who asked for medical care with symptoms like fever and coughing. Thus, people who did not seek medical attention were tested on rare occasions \cite{ilmessaggero}. Hence, for any practical means the compartment $I_S$ reflects the cases with mild to more severe symptoms. \par
In our model:  a) the number of mild to severe symptomatic cases as represented by the compartment $I_S$ does not coincide with the reported/confirmed number of infected cases as: (i) not all symptomatic cases are reported as confirmed cases, (ii) the confirmed reported total positive cases is the sum of the hospitalized and guaranteed cases, b) the number of confirmed recovered cases does not coincide with the total number of recovered cases dismissed from the hospital $R_H$; the reported cases include also a part of the confirmed cases that were isolated at home and after a second test were identified as recovered \cite{ilpost}; thus we expect that the number of $R_H$ to be less that the reported one, c) there are no data for the number of recovered from the compartment $Q$; since these cases experience mild symptoms, we assume there are no deaths from this compartment. Hence, in our model only the sum of the $H$ and $Q$ compartments reflect the number of actual confirmed cases.\par
For our analysis, and for such a short period, we assume that the total number of the population remains constant. Based on demographic data, the total population of Lombardy is $N=10m$.
The rate at which a susceptible ($S$) becomes exposed ($E$) to the virus is proportional to the density of asymptomatic infectious persons $I_A(t)$ in the total population and the symptomatic ones $I_S(t)$ the period preceding their hospitalization or home quarantine. The proportionality constant is the ``effective" disease transmission rate, say $\beta=\bar{c}p$, where $\bar{c}$ is the average number of contacts per day and $p$ is the probability of infection upon a contact between a susceptible and an infected.\par

The assumption of the model is that only a fraction, $\epsilon$, of the actual number of exposed cases $E(t)$ will experience mild to more severe symptoms and ask for medical treatment; this fraction is represented by the compartment $I_S(t)$. 
Thus, our discrete mean field compartmental SEASQHRD model reads:
\begin{align}
 &S(t)=S(t-1)-\frac{\beta(I_A(t-1)+I_S(t-1))S(t-1)}{N-D(t)-Q(t)-H(t)} \label{eq1}\\ 
 &E(t)=E(t-1)+\frac{\beta(I_A(t-1)+I_S(t-1))S(t-1)}{N-D(t)-Q(t)-H(t)}  
-\sigma E(t-1) \label{eq2}\\ 
 &I_A(t)=I_A(t-1)+(1-\epsilon) \sigma E(t-1)-\delta_A I_A(t-1)\label{eq3}\\
 &I_S(t)=I_C(t-1)+\epsilon \sigma E(t-1)-h I_S(t-1)\label{eq4}\\
 &H(t)=H(t-1)+f_1 h I_S(t-1)- f_2 \gamma H(t-1) -(1-f_2)\delta_H H(t-1) \label{eq5}\\
 &Q(t)=Q(t-1)+(1-f_1) h I_S(t-1)-\delta_Q Q(t-1)\label{eq6}\\
 &R_A(t)=R_A(t-1)+\delta_A I_A(t-1) \label{eq7} \\
 &R_H(t)=R_H(t-1)+(1-f_2) \delta_H H(t-1) \label{eq8} \\
 &R_{Q}(t)=R_{Q}(t-1)+\delta_Q Q(t-1)\label{eq9}\\
 &D(t)=D(t-1)+f_2 \gamma H(t-1)\label{eq10}
\end{align}
The model is defined in discrete time points $t=1,2,\ldots$, with the corresponding initial condition at the very start of the pandemic ($t_0$=0=DAY-ZERO) being $S(0)=N-1$, $I_{A}(0)=1$, $E(0)=I_{C}(0)=Q(0)=H(0)=R_{A}(0)=R_{H}(0)=R_{Q}(0)=D(0)=0$.\par
The values of the epidemiological parameters were fixed in the proposed model at values reported by clinical studies as follows:
\begin{itemize}
    \item $\beta (d^{-1})$ is the ``effective" transmission rate of the disease; recently published studies have reported that this rate is for any practical means the same for both asymptomatic and symptomatic cases (see e.g. \cite{jamapediatrics.2021.2025}). The transmission rate $\beta$ cannot be obtained by clinical studies, but only by mathematical models. Here we have set $\beta=0.68$ based on a previous work for the region of Lombardy\cite{Russo2020}.
    \item $\sigma (d^{-1})$ is the average per-day ``effective" rate at which an exposed person becomes infectious. In many studies, this is set equal to the inverse of the mean incubation period (i.e. the time from exposure to the development of symptoms). A study in China \cite{LIincubdoi:10.1056/NEJMoa2001316} suggests that it may range from 2 to 14 days, with a median of 5.2 days. Another study using data from 1,099 patients with laboratory-confirmed 2019-nCoV ARD from 552 hospitals in 31 provinces/provincial municipalities of China suggested that the median incubation period is 4 days (interquartile range, 2 to 7). \textit{However, the incubation period does not generally coincide with the time from exposure to the time that someone starts to be infectious.} Regarding COVID-19, it has been suggested that an exposed person can be infectious well before the development of symptoms \cite{He_2020}. Indeed, it has been suggested that for a mean incubation period of 5.2 days infectiousness starts from 2.3 days (95\% CI, 0.8–3.0 days) before symptom onset \cite{He_2020}. In our model, as explained above, $\nicefrac{1}{\sigma}$ represents the period from exposure to the onset of the contagious period and according to the above has been set as $\sigma=\nicefrac{1}{3}$.
    \item $h (d^{-1})$ is the average per-day rate at which an infected person is hospitalized. The inverse of the rate of hospitalization $h$, i.e. $\nicefrac{1}{h}$ is the mean time from the onset of symptoms to hospitalization for those infected cases that need hospitalization. For Lombardy, the median value of this parameter has been reported to be 7 days \cite{ISSDec2020}. Thus, considering also the fact that the infectiousness periods start 2-3 days before the onset of symptoms during the incubation period as explained above, we have set $\nicefrac{1}{h}={4+3}$ and thus $h=\nicefrac{1}{7}$. This coincides with the corresponding number reported by the Istituto Superiore di Sanit\a', Italy \cite{ISSDec2020}.
    \item $\delta_A (d^{-1})$ is the average per-day ``effective" recovery rate within the group of asymptomatic cases in the total population. Here, it was set equal to $\approxeq\nicefrac{1}{6.6}$ based on studies reporting that infectiousness declines quickly within $\sim$7 days \cite{He_2020}; this coincides with the mean value of the serial interval reported for Lombardy (6.6 days) \cite{cereda2020early} but also that for China (mean serial interval of 6.3 days (95\% CI 5·2–7·6) )\cite{Bi_2020}. 
    \item $\delta_H (d^{-1})$ is the average per-day ``effective" recovery rate within the subset of hospitalized cases ($H$) that finally recover. As reported recently\cite{grosso2021decreasing}, the overall median length of stay (LoS) in hospitals, which for any practical means it coincides with the recovery period, decreased steadily from 21.4 (20.5–22.8) days in February to 5.2 (4.7–5.8) days in June. Based on the above, we have set the recovery period for the hospitalized cases as $\delta_H=\nicefrac{1}{24}$ for the period until March 8, $\delta_H=\nicefrac{1}{18}$ for the period March 9-March 20, $\delta_H=\nicefrac{1}{14}$ for the period March 21-April 10, and $\delta_H=\nicefrac{1}{10}$ for the period April 11-May 4.
    \item $\delta_Q (d^{-1})$ is the average per-day ``effective" recovery rate within the subset of quarantined cases ($Q$). In a study that is based on 55,924 laboratory-confirmed cases, the WHO-China Joint Mission has reported a median time of 2 weeks from onset to clinical recovery for mild cases \cite{Chen2020n}. Hence, based on these reports, we have set $\delta_Q=\nicefrac{1}{14}$ for the quarantined cases.
    \item $\gamma (d^{-1})$ is the average per-day ``effective" mortality rate within the subset of $H(t)$ hospitalized cases that finally die. The median time from hospitalization to death for Lombardy during the early phase of the first wave of the pandemic (until April)  has been reported to be seven days through the first wave, with an inter-quartile range (3-13) \cite{ISS,grosso2021decreasing}. We should also note that the time from hospitalization to death subsequently increased as medical personnel gained experience. Thus based on the above the time from hospitalization to death was set as $\gamma=\nicefrac{1}{5}$ for the period until March 20, $\gamma=\nicefrac{1}{7}$ for the period March 21-April 10, and $\delta_H=\nicefrac{1}{10}$ for the period April 11-May 4.
    \item $f_1$ is the ratio of  all medium to severe infected cases (represented by the $I_S$ compartment) that gets hospitalized. For the early phase of the pandemic (until March 8th), a 50\% over a total 5,830 confirmed cases reported in Lombardy were hospitalized \cite{cereda2020early}, this percentage is also confirmed by other studies (see e.g. \cite{sartor2020covid}). Regarding the next period (March 20$^{\rm th}$ to May 4), that was characterized by a saturation of the hospitals beds  it is reasonable to expect that this ratio would fall significantly in addition due to the fear of a high-risk nosocomial infection. Indeed, for the period March 20 and on this ratio was about 20\% as reported also by other studies (see e.g. \cite{sartor2020covid}). Based on these studies and through the analysis of the provided data, this ratio was set to $f_1=0.65$ for the period until March 8, $f_1=0.60$ for the period March 9 to March 20, $f_1=0.5$ for the period March 21 to April 10, and $f_1=0.2$ for the period April 11 to May 4.
    \item $f_2$ is the case fatality ratio of all hospitalized cases. Based on a study considering 38,715 hospitalized COVID-19 patients in Lombardy between 21 February and 21 April 2020, this ratio was relatively high \cite{ferroni2020survival,grosso2021decreasing}) mainly due to the saturation of ICU beds. In particular as reported in a recent published study \cite{grosso2021decreasing} for the region of Lombardy is estimated to have steadily decreased from 34.6\% (32.5–36.6\%) in February to 29.8\% (29.3–30.3\%) in March, 22.0\% (21.1–23.0\%) in April, 14.7\% (13.3–15.9\%) in May and 7.6\% (6.3–10.6\%) in June. Based on these studies, the fatality ratio was set $f_2=0.27$ for the period February 24-March 8, $f_2=0.23$ for the period March 9-April 10, and $f_2=0.20$ for the period April 11 -May 4.
    \item $\epsilon (d^{-1})$ is the fraction of the exposed cases that enter compartment $I_S$, i.e. those that will experience mild to severe symptoms and ask for medical care. For the case of Italy and in particular that of Lombardy during the first wave of the pandemic; the ratio of the confirmed to unreported cases was estimated to be around 0.07 (interquartile range: 0.023 to 0.1), \cite{Russo2020}. In our model, we have set $\epsilon=0.12$ for the period until March 8,  $\epsilon=0.10$ for the period March 9 to March 20, $\epsilon=0.09$ for the period March 21 to April 10, and $\epsilon=0.05$ for the period April 11 to May 4.
\end{itemize}
%As new daily confirmed cases of recovered and dead appear with a time delay (which is generally unknown but an estimate can be obtained by clinical studies) with respect to the corresponding infected cases, the above per-day rates are not the actual ones; thus, they are denoted as ``effective" daily rates.
Regarding DAY-ZERO in Lombardy, i.e. the initial date of the introduction of the disease in the region, that was set to January 15, as estimated based on a methodological framework based on genetic algorithms together with the level of asymptomatic cases in the total population \cite{Russo2020}. Another study \cite{Zehender} based on genomic and phylogenetic data analysis, reports the same time period, i.e. that between the second half of January and early February, 2020, as the time when the novel coronavirus SARS-CoV-2 entered northern Italy.

\subsection{Estimation of the basic and effective reproduction numbers}\label{subsubsec22}

Note, that there are three infected compartments, namely $E,I_A,I_S$ that determine the fate of the pandemic. Thus, considering the corresponding equations given by Eq.\eqref{eq2},\eqref{eq3},\eqref{eq4}, and that at the very first days of the epidemic $S\approx N$, the Jacobian of the system as evaluated at the disease-free state is:
\begin{align}
\boldsymbol{J} & \triangleq \frac{\partial(E(t),I_A(t),I_S(t)) }{\partial(E(t-1),I_A(t-1),I_S(t-1)}\nonumber\\
& =
\begin{bmatrix}
1-\sigma & \beta & \beta\\
(1-\epsilon)\sigma & 1-\delta_A & 0\\
\epsilon\sigma & 0 & 1-h
\end{bmatrix}
=
\begin{bmatrix}
1 & 0 & 0\\
0 & 1 & 0 \\
0 & 0 & 1
\end{bmatrix}
+
\begin{bmatrix}
-\sigma & \beta & \beta\\
(1-\epsilon)\sigma & -\delta_A & 0\\
\epsilon\sigma & 0 & -h
\end{bmatrix}\label{eq11}
\end{align}
The eigenvalues (that is the roots of the characteristic polynomial of the Jacobian matrix) dictate if the disease-free equilibrium is stable or not, that is if an emerging infectious disease can spread in the population. In particular, the disease-free state is stable, meaning that an infectious disease will not result in an outbreak, if and only if all the norms of the eigenvalues of the Jacobian $\boldsymbol{J}$ of the discrete time system  are bounded by one. Jury's stability criterion \cite{Jury_1976} (the analogue of Routh-Hurwitz criterion for discrete-time systems) can be used to determine the stability of the linearized discrete time system by analysis of the coefficients of its characteristic polynomial. The characteristic polynomial of the Jacobian matrix is:
\begin{equation}
    F(z)=a_3z^3+a_2z^2+a_1z+a_0
\end{equation}
\noindent where
\begin{align}
& a_3=1\nonumber\\
& a_2=\delta_A + h + \sigma - 3 \nonumber\\
& a_1=\delta_A h - 2h - 2\sigma - 2\delta_A - \beta \sigma + \delta_A \sigma + h \sigma + 3 \nonumber\\
& a_0= \delta_A + h + \sigma - \delta_A h + \beta \sigma - \delta \sigma - h \sigma - \beta h \sigma + \delta_A h \sigma + \beta \epsilon h \sigma - \beta \delta_A \epsilon \sigma - 1
\end{align}
\noindent The necessary conditions for stability read:
\begin{align}
F(1)>0  \label{eq15}\\
(-1)^3 F(-1)>0 \label{eq16}
\end{align}
Thus, the sufficient conditions for stability are given by the following two inequalities:
\begin{align}
\vert a_0 \vert<a_3  \label{eq17}\\
 \vert b_0 \vert> \vert b_2 \vert, \label{eq18}
\end{align}
where, 
\begin{align}
b_0=\left|\begin{matrix}
a_0 & a_3\\
a_3& a_0\\
\end{matrix}\right|,%vmatrix doesn't work with IJRNC
b_2=\left|\begin{matrix}
a_0 & a_1\\
a_3& a_2\\
\end{matrix}\right|\label{eq19}
\end{align}
It can be shown, that the second necessary condition of \eqref{eq16} and the first sufficient condition of \eqref{eq17} are always satisfied for the range of values of the epidemiological parameters considered here.
The first inequality \eqref{eq15} results in the necessary condition:
\begin{equation}
\beta  (\frac{1-\epsilon}{\delta_A}+\frac{\epsilon}{h}) <1 \label{eq21}
\end{equation}
It can be also shown, that for the range of the parameters considered here, the second sufficient condition of \eqref{eq18} is satisfied if the necessary condition of \eqref{eq21} is satisfied. Thus, the necessary condition of \eqref{eq21} is also a sufficient condition for stability. Hence, the disease-free state is stable, if and only if, condition of \eqref{eq21} is satisfied.

Note that in this necessary and sufficient condition of \eqref{eq21}, the first term reflects the contribution of compartment $I_A$ 
and the second term the contribution of compartment $I_S$ in the spread of the disease. 
Thus, the above expression reflects the basic reproduction number $R_0$ which is qualitatively defined by $R_0=\nicefrac{\beta}{\textit{infection time}}$. Hence, our model results in the following expression for the basic reproduction number:
\begin{equation}
R_0=\beta  (\frac{1-\epsilon}{\delta_A}+\frac{\epsilon}{h})  \label{eq22}
\end{equation}
Note that for $\epsilon=0$, the above expression simplifies to $R_0$ for the simple SIR model.
Based on the values that we have accounted in our model, the above expression results to $R_0=4.38$. This estimation is close to the one reported in Russo et al. \cite{Russo2020} ($R_0=4.56$) and greater than the ones reported in other studies. For example, Allieta et al. \cite{allieta2021covid} reported $R_0=3.88$, while in another study regarding Italy, D’Arienzo and Coniglio \cite{d2020assessment} used a SIR model to fit the reported data in nine Italian cities and found that $R_0$ ranged from 2.43 to 3.10. Finally, \cite{pillonetto2021tracking} reported an $R_0$ around 3.4. However, we note that the above studies that provide estimates of $R_o$ are based  on the reported cases which, may lead to an underestimation of the basic reproduction number (see also the discussion in \cite{Cori_2013}), while our model takes into account the compartment of asymptomatic cases which transmit the disease. \par

Regarding the estimation of the effective reproduction number $R_e$, this corresponds to the average number of secondary infections from a single infectious individual during an epidemic. For its calculation, we derive the  next generation matrix $G$ with elements $g_{ij}$ formed by the average number of secondary infections of type $i$ from an infected individual of type $j$, given by (see also \cite{Russo2020}):
\begin{equation}
    \boldsymbol{G}=\nabla \boldsymbol{F} \cdot \nabla (\boldsymbol{V})^{-1}.
\end{equation}
$\boldsymbol{F}$ is the vector containing  the transmission rates, and $\boldsymbol{V}$ is the vector containing the transition rates between the infected compartments. The effective reproduction number $R_e$ is the spectral radius, i.e. the dominant eigenvalue of $G$.
In our model we have:
\begin{equation}\label{e:Reff}
    R_e(t)=\beta  (\frac{1-\epsilon}{\delta_A}+\frac{\epsilon}{h}) \frac{S(t-1)}{N-D(t-1)-Q(t-1)-H(t-1)}.
\end{equation}
%\begin{align}
%\boldsymbol{F}=\begin{bmatrix}
%\beta \frac{S(t-1)}{N-D(t)-Q(t)-H(t)} (I_A(t-1)+I_S(t-1))\\
%0
%\end{bmatrix}, \boldsymbol{V}=\begin{bmatrix}
%\sigma E(t-1)\\-(1-\epsilon) \sigma E(t-1)+\delta_A I_A(t-1)
%\end{bmatrix}.
%\end{align}

\subsection{Google Mobility Report Analysis}\label{s:Mobility}
In the spring of 2020, many nations enacted policies restricting resident movement in hopes of slowing the spread of COVID-19. Although the individual compliance to these policies may vary, there is a strong link between government mobility policies and the subsequent mobility behavior in a population \cite{Armstrong2020,buonomo2020effects}. Furthermore, these mobility patterns are correlated with changes in COVID-19 cases \cite{badr_association_2020}.

There are several dates that demarcate major changes in Lombardy mobility policies. On March 4, 2020 all schools and universities were closed. On March 8, 2020, the first major lockdown was introduced, with the encouragement of working from home and closure of recreational venues. On March 21, 2020, measures were furthered; a stay-at-home order was issued along with the banning of pedestrian activity. These measures continued until May 4, when policies were gradually relaxed; non-essential industrial activities reopened, and other pedestrian activities was allowed. These key dates can be used to divide the initial wave of COVID-19 into four approximately time periods. The first period- February 15 to March 8- starts with the availability of mobility data in Lombardy and ends with the introduction of lockdown. The second period stretches from March 9 to March 20: the period of initial lockdown. The third period represents the first phase of stricter measures, ranging from March 21 to April 19, while the fourth and final period begins April 20 and continues until policy relaxation on May 4. 

To describe mobility patterns in Lombardy, we used data provided by the Google Community Mobility Reports \cite{Google}. The reports show movement trends within a population for a given region. The mobility data is organized into six categories (i.e. activities): retail \& recreation (RR), grocery/pharmacy (G), parks (P), transit (T), workplaces (W), and residential (R). We note that within this period schools \& universities (SU) were closed. The mobility values are expressed as percent deviation from the baseline value; the baseline was set as the median value from January 3 to February 6, 2020. All activities except for the residential one determine the percent deviation from the change in total visitors. The residential category is determined by the change in time spent in places of residence. 

The relative impact of each mobility activity on disease transmission was analyzed to inform the parameter $\beta$ in the model, since the amount of personal contacts made in a given location can be used to estimate the likelihood of disease transmission in that location. Based on a survey by Mossong et. al. \cite{mossong_social_2008}, home, work, and leisure activities each account for about 20\% of all personal contacts reported, transport accounts for the least amount of contacts at about 5\% and a 20\% is accounted to multiple activities. Despite travel being responsible for the least amount of physical and verbal contacts, diminished travel has still been strongly correlated with a decrease in the reproduction number \cite{vollmer_report_2020}. Therefore, based on this information, each reported mobility category likely plays a significant role in the transmission of COVID-19.
%\begin{center}
\begin{table}[ht]%
\centering
\caption{Analysis of Google mobility report's data in Lombardy during the first wave of the pandemic (expected value of maximum mobility restriction and associated standard deviation of each activity)}
\label{t:Google_mobility}%Original at end of tex
\begin{tabular*}{500pt}{@{\extracolsep\fill}ccccccc@{\extracolsep\fill}}
\toprule
 \textbf{Retail\& Recreation (RR)}  & \textbf{Grocery} & \textbf{Parks} & \textbf{Transit} & \textbf{Workplaces} & \textbf{Residential}\\
\midrule
-91.4\% (6.5) &-58.6\% (11) &-84.3\% (23.3)  &-86.7\% (6.4)  &-74.6\% (11.3) &34\% (4.6)\\
\bottomrule
\end{tabular*}
\end{table}
%\end{center}
To account for the role of each activity in the transmission of COVID-19, an aggregate value was found by appropriately weighting the six categories for a given day. Thus, the transmission rate is impacted by the change in mobility such that $\beta=(1-u)(1-\theta_C)\beta_0$, with $u$ a measure of the reduction in mobility, i.e., social distancing. In this formulation, when $u = 0$, there are no mobility restrictions, and when $u = 1$, there is assumed to be complete lack of mobility and social interaction in the population beyond the residential.
The caution that people exhibit in their daily interactions is captured by the term $\theta_C$, with $\theta_C=0$ being no caution. Following the first curtails in mobility on March 8$^{\rm th}$, the caution was quantified to the value $\theta_C=0.2$. This is in agreement with what has been reported also in other studies about the effect of the cautiousness in the reduction of the disease transmission rate \cite{Caley_2007}.
%For example, if the Google data demonstrates a 80\% decrease in mobility, $u$ would be assumed to be 0.8. For each of the four time periods, the data obeys an approximately Gaussian distribution. 
Analyzing the available data, the effect of policy changes for different aspects of life in society can thus be quantified as 
\begin{equation}\label{u_def}
    u(t)=0.2\,RR(t)+0.05\,G(t)+0.05\,P(t)+0.05\,T(t)+0.25\,W(t)+0.2\,SU(t)+0.2\,R(t) 
\end{equation}
Obviously, the mobility activity in the category ``Residential'' is affected by the other activities, including the transit which includes also traveling to other regions for vacation or work. In order to model and assess this dependence, a linear regression model of the following form was fit:
\begin{equation}
    R(t)=c_0+c_1\,RR(t)+c_2\,G(t)+c_3\,P(t)+c_4\,T(t)+c_5\,W(t)+c_6\,SU(t)
    \label{fithome}
\end{equation}
The parameters of the model were estimated using the available mobility data from Feb 15, to May 4 2020. In Table~\ref{t:fit}, we report the values of the parameters and their 95\% confidence intervals together with the fit statistics, namely the standard error (SE) of the coefficients, the t-statistic to test the null hypothesis that the corresponding coefficient is zero against the alternative that it is different from zero, given the other predictors in the model, the corresponding p-value for the t-statistic of the hypothesis test that the corresponding coefficient is equal to zero or not. The Root Mean Squared Error was 0.0194, the R-squared 0.979,  and the adjusted R-squared was 0.978.
%\begin{center}
\begin{table}[ht]%
\centering
\caption{Parameter estimates and statistics of the linear regression model fit of the model~(\ref{t:fit})}
    \label{t:fit}
\begin{tabular*}{500pt}{@{\extracolsep\fill}lccccc@{\extracolsep\fill}}
\toprule
 & \textbf{Estimate} & \textbf{SE}  & \textbf{$t_{stat}$} & \textbf{$p_{value}$}\\
\midrule
$c_0$ (intercept) &-0.0020 (-0.0160, 0.0120)&0.0070&-0.2866&0.7752\\
$c_1$ (RR) &0.0809 (0.0075,0.1545)&0.0369&2.1962&0.0312\\
$c_2$ (G) &0.1296 (0.0945,0.1647)&0.0176&7.3648&2.1823e-10\\
$c_3$ (P) &-0.0495 (-0.0795,-0.0194)&0.0151&-3.2784&0.0016\\
$c_4$ (T) &0.0655(-0.0151,0.1461)&0.0404&1.6205&0.1094\\
$c_5$ (W) &-0.6634(-0.7319,-0.5949)& 0.0344&-19.3091&2.1021e-30\\
$c_6$ (S) &-0.0223(-0.0430,-0.0017)&0.0104&-2.1558&0.0344\\
\bottomrule
\end{tabular*}
\end{table}
%\end{center}
%\begin{figure}[htbp]
%    \centering
%    \includegraphics[width=0.7\textwidth]{./Figures/mobility_boxplot.eps}
%    \caption{Distribution of aggregate mobility in Lombardy for the four periods. The rectangles present the range of the observed aggregate mobility with the internal line representing the median. Blue stars denote the average mobility, while the dashed lines denote the 95\% confidence range. Red crosses in the third period represent two outliers that were disregarded during the analysis.}
%    \label{f:GMD}
%\end{figure}
Thus, at a statistical significance threshold of $\alpha=0.05$, the variables RR,G,P,W and S are significant for describing the residential mobility, while the transit activity is not statistically significant. Note also that by the statistics, the hypothesis that $c_0$ is different than zero cannot be rejected at the significance level of 0.05, highlighting that the residential activity is purely dependent on the rest of the activities and can be discounted in the following section.
%\begin{figure}
%    \centering
%         (a)\includegraphics[width=0.3\textheight]{./Figures%/Distancing_OL.eps}\vspace{1em}(b)\includegraphics[width=0.3\textheight]{./Figures/hist_u.eps}
%    \caption{Social Distancing as reflected from Google Mobility Reports for the period February 24-May 4. (A) Profile in time with red solid line presenting the average value and the cyan band presenting the variance. (B) Histogram of the variance in social distancing; a standard deviation of 2.8\% is identified}
%    \label{f:GM_S1}
%\end{figure}
%Combining Eq.~\eqref{u_def} \& Eq.~\eqref{fithome},
%the resulting cumulative effect of activity curtails on the social distancing variable is presented in Fig.~\ref{f:GM_S1}. 
Employing the statistical analysis results of Table~\ref{t:Google_mobility}, the resulting standard deviation for social distancing, $u(t)$, was identified to be $\sigma_u=0.0282$ and furthermore it could be considered constant during the period Feb. 24$^{\rm th}$ to May 4$^{ \rm th}$. Note that during this period there were specific governmental policies assigning mobility curtails to each activity. From this analysis the adherence by Lombardy to the policies can be estimated, which must be accounted for.

\begin{figure}[p]
    \centering
    \subfigure[]{
         \includegraphics[width=0.29\textheight]{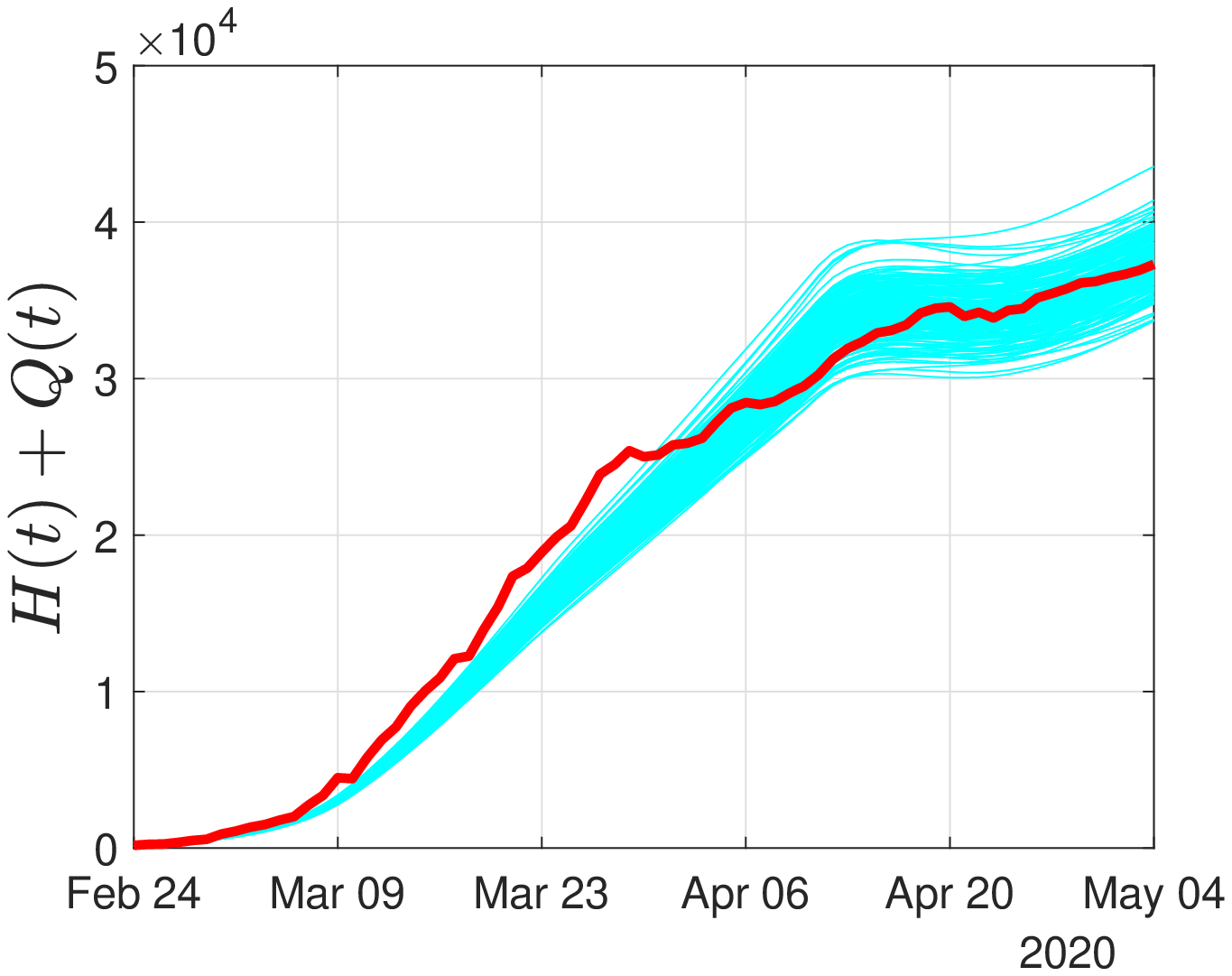}
    }
    \subfigure[]{
         \includegraphics[width=0.31\textheight]{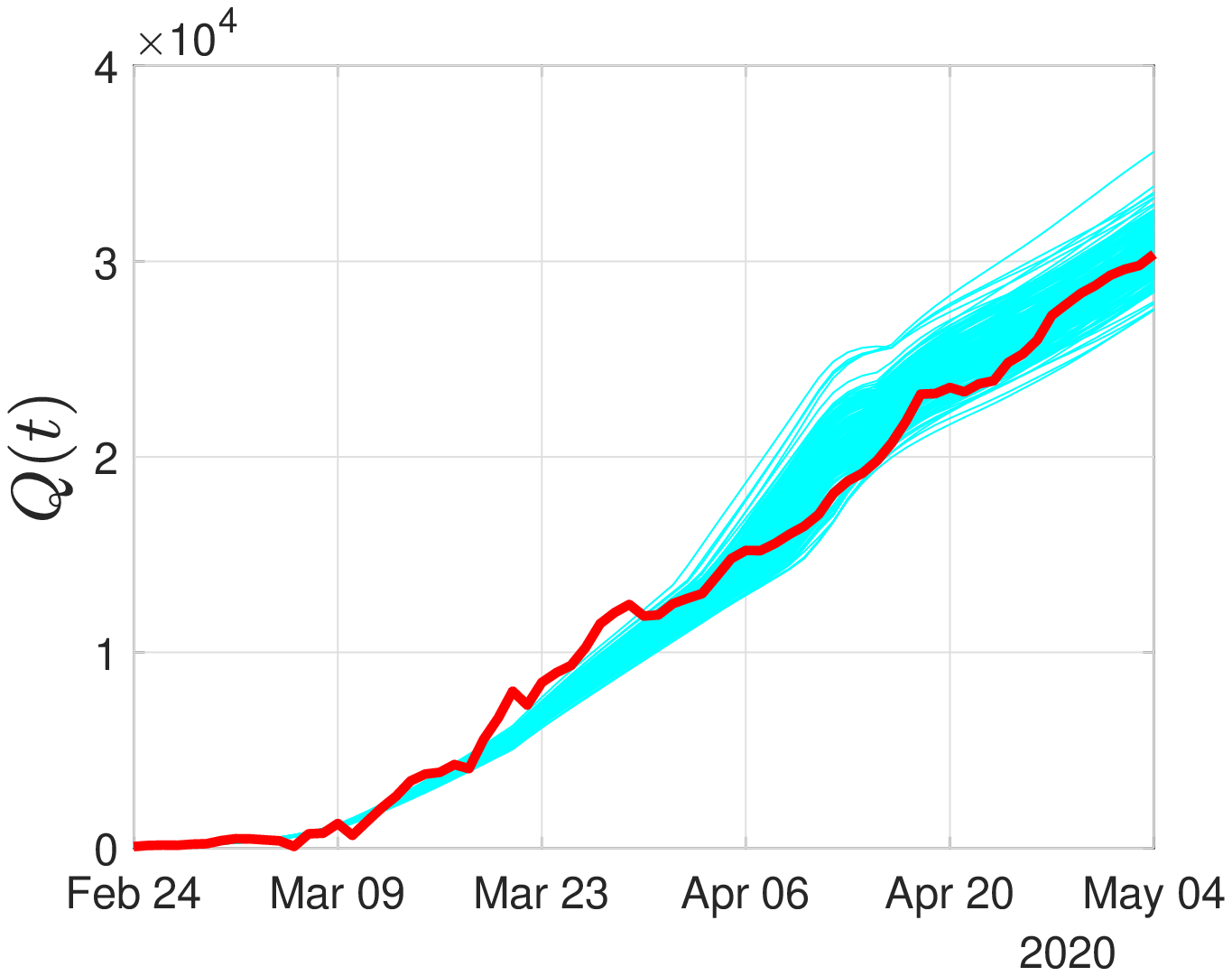}
    }
        \subfigure[]{
         \includegraphics[width=0.31\textheight]{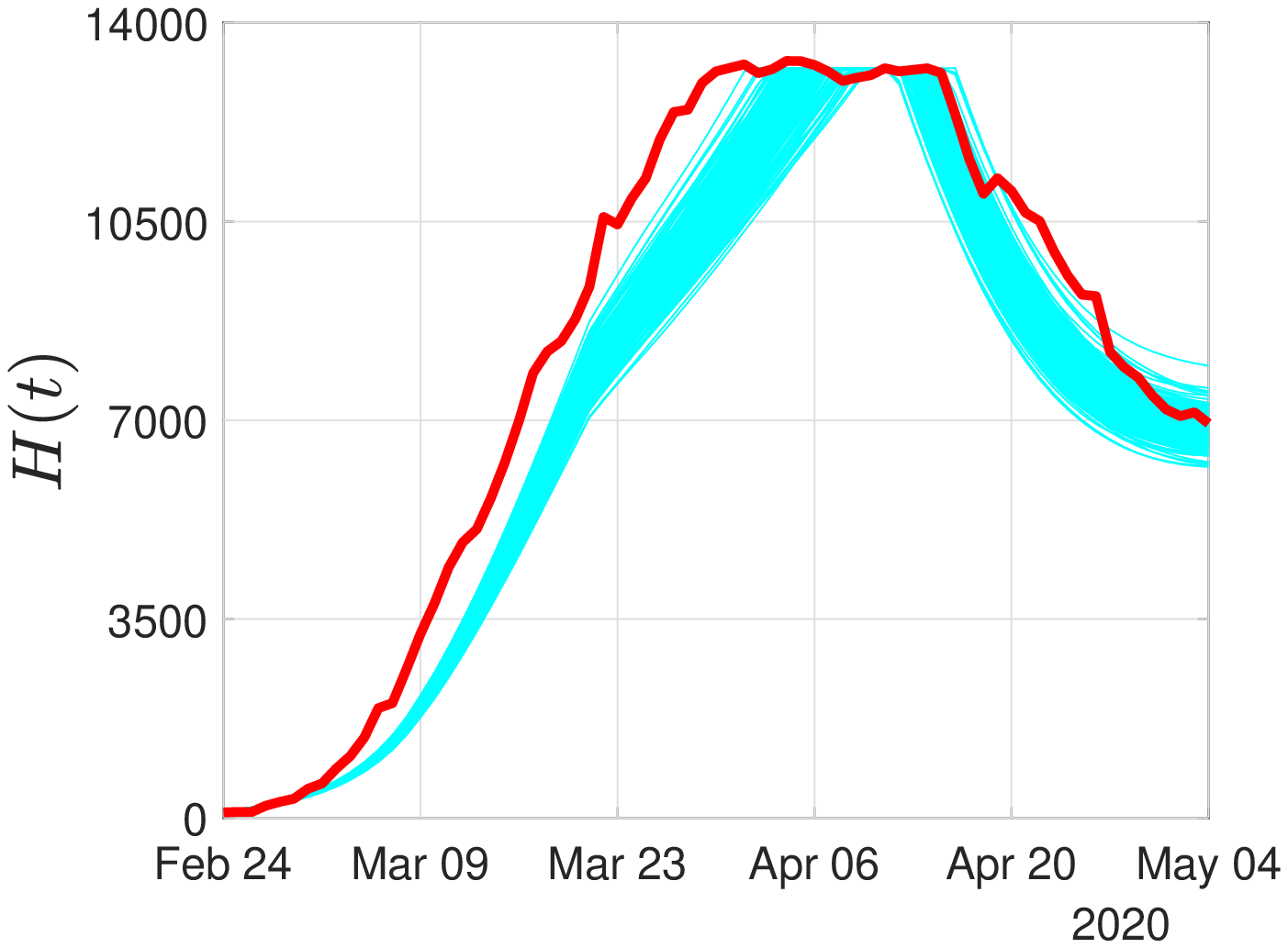}
    }
    \subfigure[]{
         \includegraphics[width=0.31\textheight]{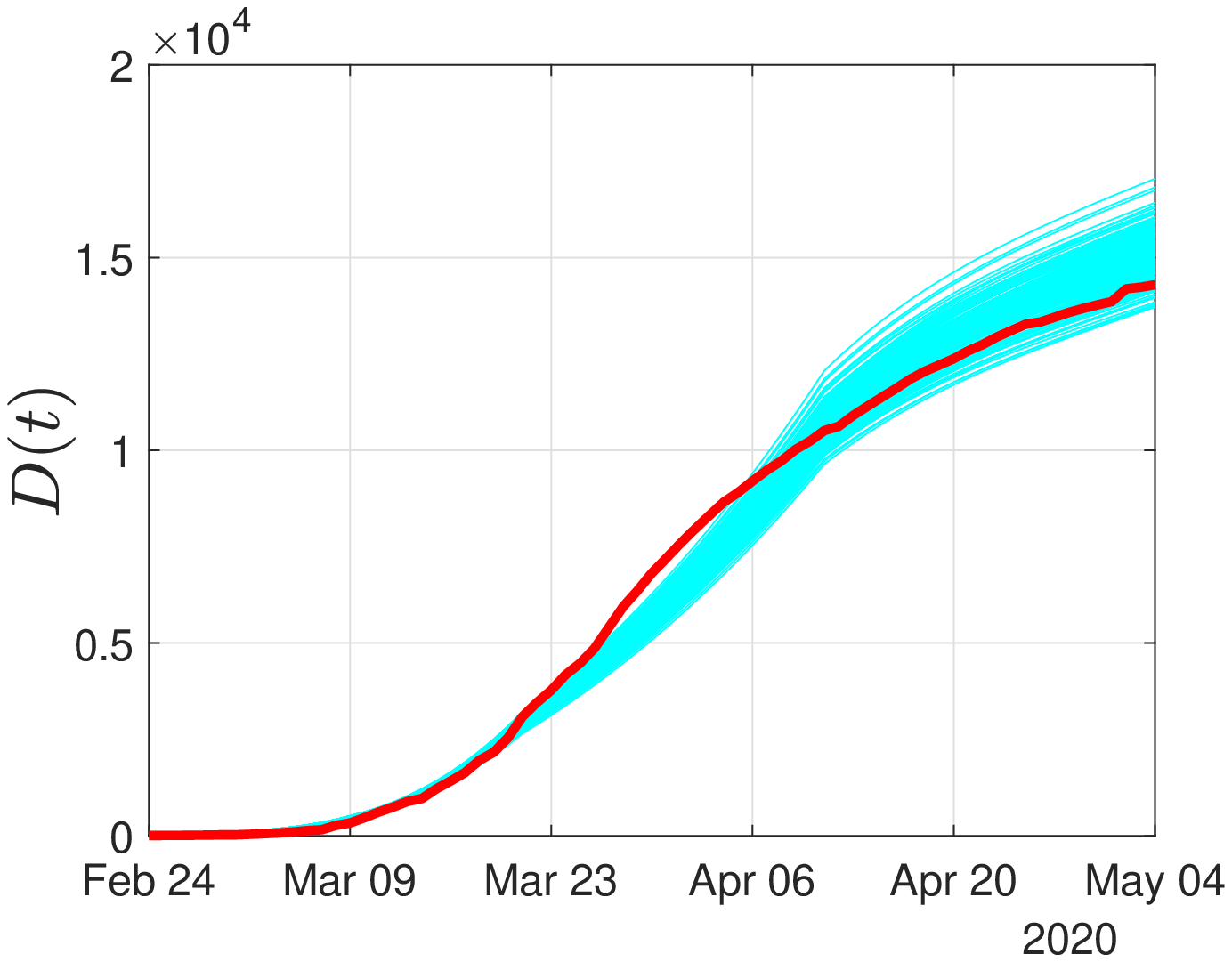}
    }
    \subfigure[]{
         \includegraphics[width=0.31\textheight]{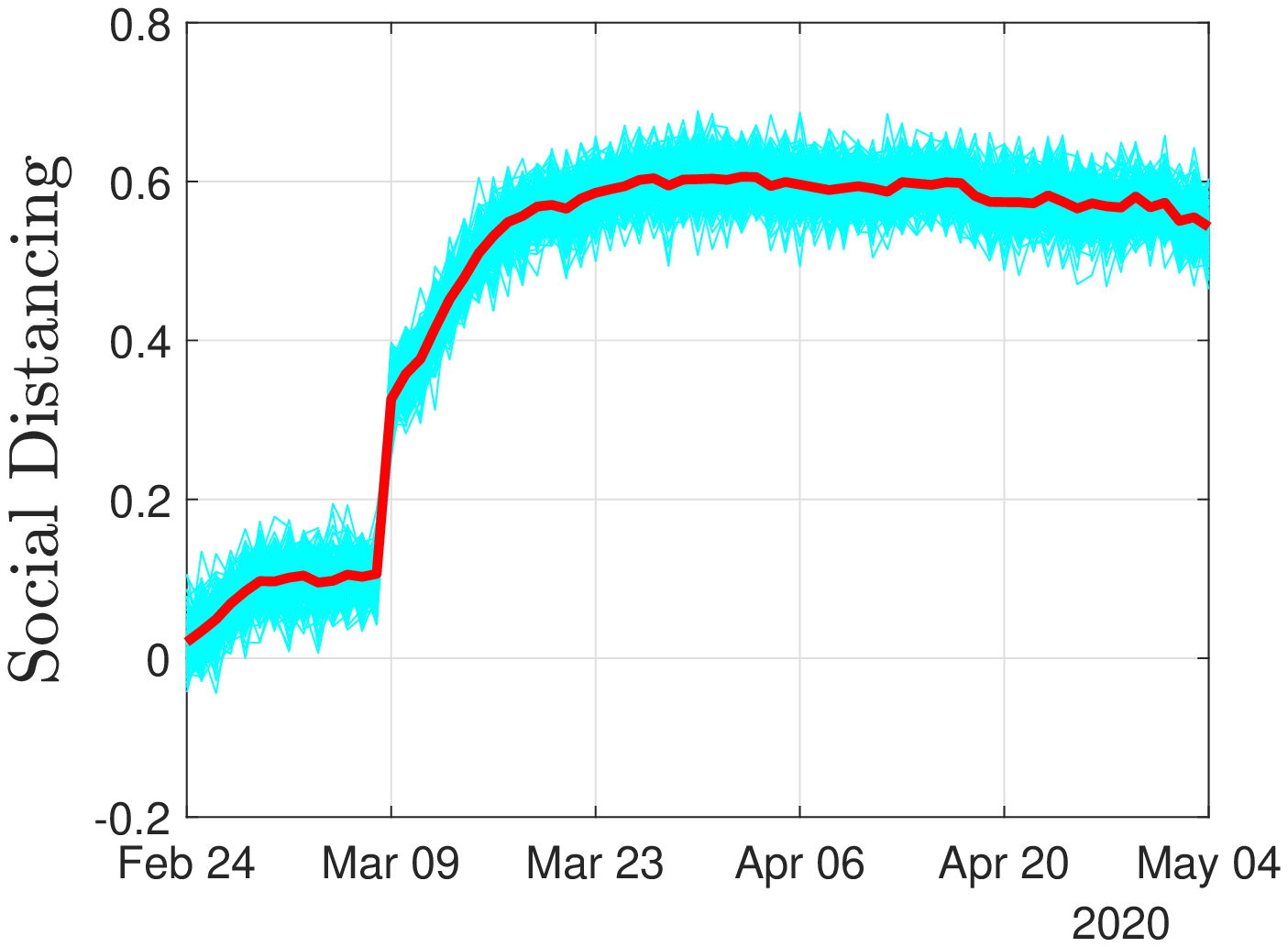}
    }
    \subfigure[]{
         \includegraphics[width=0.31\textheight]{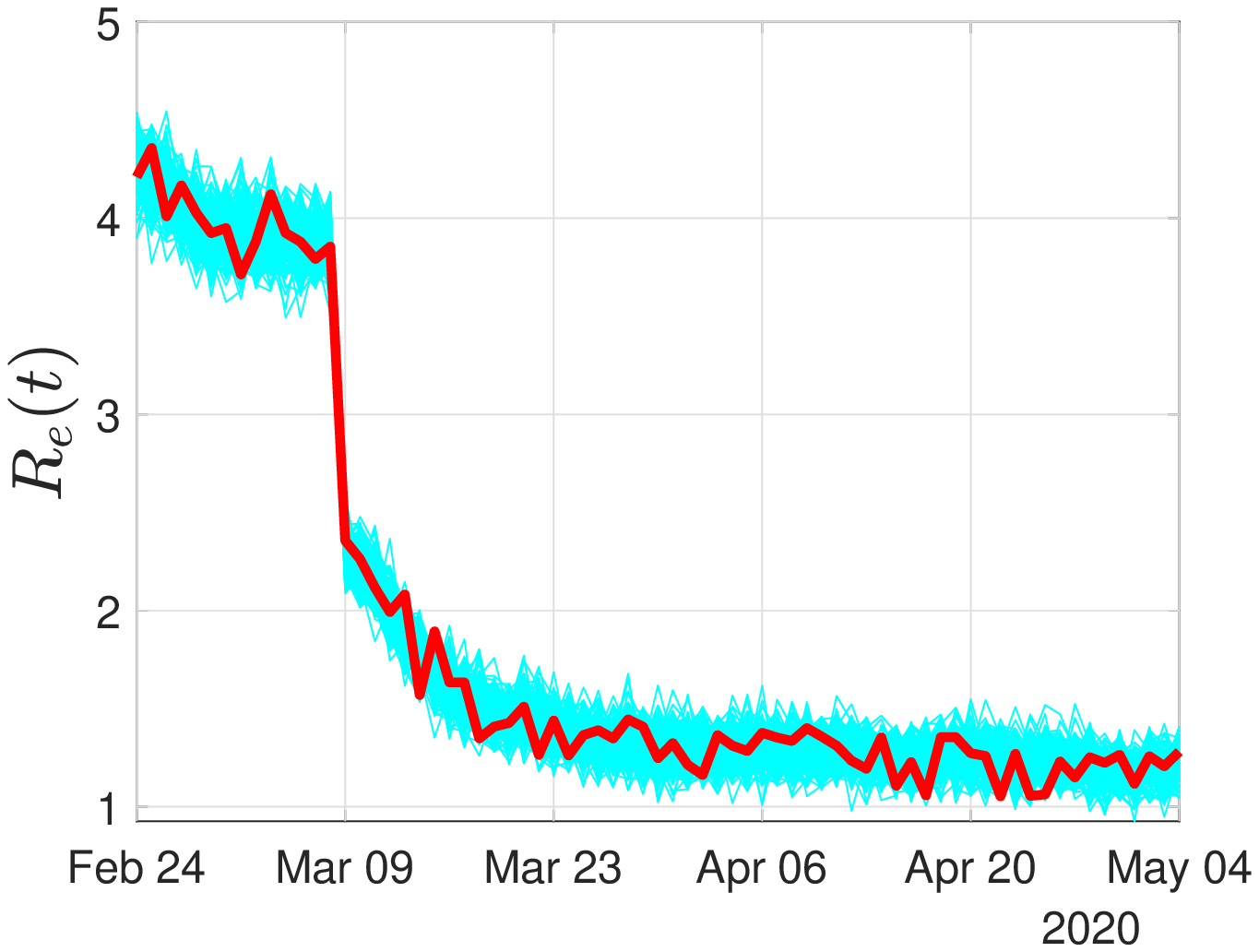}
    }
\caption{Simulation results using the SEASQHRD model considering the social distancing as obtained from Google mobility reports (see Table \ref{t:Google_mobility}) for the period February 24-May 4. (a) Reported Infected cases (in our model this corresponds to the sum of $H(t)$ and $Q(t)$). (b) Hospitalized. (c) Quarantined. (d) Dead. (e) Social distancing as reflected from Google Mobility Reports. (f) Effective Reproduction number $R_e$. Red solid lines in (a), (b), (c), (d) correspond to the actual confirmed cases. Red solid lines in (d),(e) present average values. Cyan colors present individual trajectories of 200 adherence scenaria.}\label{f:GM_1}
\end{figure}

Employing the SEASQHRD model with the identified parameters and considering the variance in social distancing due to adherence (Table 1), we simulated the progression of the disease under the social distancing policy set by the government for 200 adherence scenarios. Figure~\ref{f:GM_1} presents the simulations and the reported data for Lombardy.
As the results show, the predictions of the model fits fair the actual reported data. We note that the lower limit of $R_e$ on May 4 is close enough to 1. As discussed above, the difference between our estimates with the reported $R_e$ in other studies (see e.g. \cite{pillonetto2021tracking} who reported a value of $R_e$ around 0.6 just before the end of the lockdown on May 4) is due to the fact our model takes into account also the asymptotic cases assuming that they also transmit the disease.

\subsection{Probabilistic Model Predictive Control}\label{subsubsec24}
The developed SEASQHRD model of \S\ref{s:model} was employed to propose government policy as the pandemic progresses. A systematic mechanism for the proposed policy hinges on the solution of an optimization problem that balances the impact on quality of life due to prolonged social distancing on the one hand and the impact of the pandemic on life expectancy on the other. Considerations include the question of societal adherence to government advisories as well as other unaccounted factors which affect the outcome of predictions. We thus assume that government policies need to be continually evaluated and adapted to reality. 
A system's approach towards this end involves the design of a probabilistic model predictive control (PMPC) structure. The design formulates the policy question as an optimization problem with a finite number of policy decisions being made and their effect being predicted by the SEASQHRD simulator for a finite time horizon forward, called the prediction horizon, $T_P$, given the current state of the epidemic. Each policy decision is effective for a finite period, called the Action Horizon, $T_A$. The number of decisions, $N_d$, multiplied by $T_A$ is called the control horizon, $T_C\triangleq N_d\, T_A$. Once a policy has been identified, the first policy decision is enacted and by the end of period $T_A$ the state of the epidemic is measured. A new optimization problem for the period $T_A$ to $T_A+T_P$ is then formulated and solved. This recursion constitutes the MPC, where one of the advantages is that it can adapt to uncertainty in the epidemic's behavior. 
Note, that the state of the compartments is a prerequisite for the SEASQHRD simulator, however only three compartments are directly measured, leaving seven compartment states unknown. 
In the previous section the analysis of the
activities resulted in a distribution of values which impacted the epidemic's evolution. Recognizing the uncertainty in the announced policy due to adherence and the limitations of the derived model, the underlying optimization problem is defined as a probabilistic one
\begin{align}
&    \alpha^*(t)=\arg\min_{\delta \alpha\in\mathbb{R}^{6\times N_d}}
\sum_{t=\text{Day 1}}^\text{Day $T_P$} \omega_I {\rm \bf E}\left[\left(\frac{H(t;u,\theta_A)}{H(t;u\equiv 0,\theta_A)}\right)^2\right] 
+ \omega_R {\rm \bf E}\left[ R_{eff}(t;u,\theta_A)^2\right]
+ \alpha(t)^T\Omega_A\alpha(t)\label{OP1:J}\\
&    \text{subject to}\nonumber\\
&    \alpha(t)=\sum_{j=1}^{N_d} {\bf H}(t-(j-1)\,T_A)\delta \alpha_{j}\label{OP1:u}\\
&    {\bf 0}\le \alpha(t)\le [0.91\; 0.59\; 0.85\; 0.87\; 0.75\; 1]^T \label{OP1:0}\\
&    \delta\alpha_j\le \delta\alpha_c \label{OP1:0a}\\
&    u(t)=[0.216\;0.076\;0.04\;0.063\;0.117\;0.196]\,\alpha(t) \label{OP1:_1}\\
&    \beta(t)=\beta_0(1-\theta_C(t))(1-u(t)+\theta_A(t))\\
&    {\rm \bf P}[H(t; u,\theta_A) \le Beds]\le P^c_{bed}, \quad t=\{1,...,T_p\}\label{OP1:1}\\
&    {\rm \bf P}[R_{eff}(t; u,\theta_A) \le R_{e,c}]\le P^c_{rep}, \quad t=\{1,...,T_p\}\label{OP1:2}\\
&    X(t;u,\theta_A)= {\rm SEASQHRD}(X(t-1;u,\theta_A),u(t),\theta_A),\;\; X(0; u)=X_0,\; \; t=\{1,...,T_p\}\label{OP1:3}
%&    u(t)=0.216\,RR(t)+0.076\,G(t)+0.04\,P(t)+0.063\,T(t)+0.117\,W(t)+0.196\,S(t) \label{OP1:_1}\\
%    
\end{align}
where $\alpha(t)$ is the governmental distancing policy, consisting of curtails on activities 
$\alpha(t)= [RR(t)\;G(t)\;P(t)\;T(t)\;W(t)\;SU(t)]^T$; each activity in $\alpha$ is normalized to map $(0,T_P] \mapsto [0,1]$, where $0$ implies no restrictions on the corresponding activity and $1$ implies complete social curtail on it. Function ${\bf H}(\cdot)$ denotes the standard Heaviside function. Equation \eqref{OP1:u} denotes that during action period $j$ policies in $\alpha_j$ are active. Following standard MPC practice the final policies $\alpha_{N_d}$ are maintained till the end of the prediction horizon; one reason for this practice is that a final time penalty term can then be neglected. The state of infection, $X$, in \eqref{OP1:3} is defined as $X=[S\,E\,I_A\,I_S\,H\,Q\,R_A\, R_H\,R_Q\,D]$. Term ${\bf 0}$ denotes the zero vector of appropriate size. In the specific PMPC formulation a number of epidemic progression scenaria (defined as $N_{cases}$) is investigated; the adherence of Lombardy's population to the government curtails is the uncertain parameter which is sampled based on the analysis results of Table \ref{t:Google_mobility}. The trajectory statistics are then computed to obtain the distributions used in the formulation.

The cost function of \eqref{OP1:J} balances the social cost of the epidemic to the social distancing, via the relative weights $\omega_I$ and $\Omega_A$. Function ${\rm\bf E}[\cdot]$ denotes the expectation of the distribution for argument $\cdot$, while ${\rm\bf P}[\cdot]$ denotes the probability argument $\cdot$ is true. The cost function is constructed to be convex to facilitate convergence to an optimum. 
Constraint \eqref{OP1:0} reflects the need for mobility of essential workers and for essential activities to take place in order for society to keep functioning. These upper bounds where computed from the Google mobility report data (presented in Table~\ref{t:Google_mobility}), by estimating the maximal expected curtail on mobility during the period Feb-24 to May 4, during which society was in essence completely closed. 

The stress on society due to mobility restrictions can be significant. Constraint \eqref{OP1:0a} is introduced to prevent a rapid increase in the restriction of each activity. In the original formulation, we impose a constraint of 25\% increase in each activity restriction with the exception of schools \& universities ($SU$) which is left unconstrained (i.e., $\delta\alpha_c=[0.25\;0.25\;0.25\;0.25\;0.25\;1.0]$). Note that there is no constraint imposed for lifting the restrictions.
%\begin{figure}[htbp]
%    \centering
%         (a)\includegraphics[width=0.3\textheight]{./Figures/Reff_OL.eps}
%    \caption{Effect of Social Distancing as reflected from Google Mobility Reports for the period February 24-May 4 on the effective reproduction number; the red solid line presents the expected value and the cyan band presents the variance.}
%    \label{f:Reff_OL}
%\end{figure}
The societal effect of the desired policy is also captured by the use of two constraints, \eqref{OP1:1} limiting the overflow of patients to hospitals and \eqref{OP1:2} limiting the effective reproduction number, defined in \eqref{e:Reff} (with the basic reproduction number, $R_0$ calculated from \eqref{eq22}) so that new infection cases dynamic becomes a contraction. In Figure~\ref{f:GM_1}(f), we observe that during the period February 24-May 4, $R_e$ remains above one (even though it is greatly reduced) which hints at the increasing rate of asymptomatic cases that transmit the disease during the same period.
The use of probabilities in constraints \eqref{OP1:1} \& \eqref{OP1:2}
captures the degree of risk aversion in the desired policy, by imposing a limit to the risk of violating them via tunable parameters $P^c_{bed}\in[0,1]$ and $P^c_{rep}\in[0,1]$. A value of $0.5$ converts this constraint to an expectation constraint. For the specific formulation, as the values of $P^c_{bed}$ and $P^c_{rep}$ increase the PMPC formulation converges to the robust MPC formulation\cite{MHASKAR2005209}.
\section{Simulation Results}

The model predictive control simulation begins on February 24$^{\rm th}$, three days after the first case was reported in Lombardy. 
Within this 70 day period until May 4, when the restrictions were relaxed, the horizon lengths determine the duration and frequency of policy investigation and implementation. The prediction horizon, or how far into the future the simulation investigates at a given time point, was examined at lengths of 1 to 4 weeks. As the prediction horizon increases from 1 to 4 weeks, the action during the $1^{\rm st}$ period converges to a set value. A prediction horizon of 3 weeks was chosen to balance this convergence and the uncertainty in model predictions. The control horizon, which determines the number of policies during a given prediction, was chosen as 2 weeks to balance problem complexity and solution convergence. Figure \ref{f:MPC_S1} of the Supporting Information showcases the effect of the different horizon choices. The action horizon, that is the actual duration of a given policy, was chosen to be 1 week. This choice of horizon allows for frequent policy changes if necessitated by the model predictions, while not being so frequent as to be unrealistic from a policy implementation standpoint.

\begin{figure}[p]
    \centering
    \subfigure[]{
         \includegraphics[width=0.31\textheight]{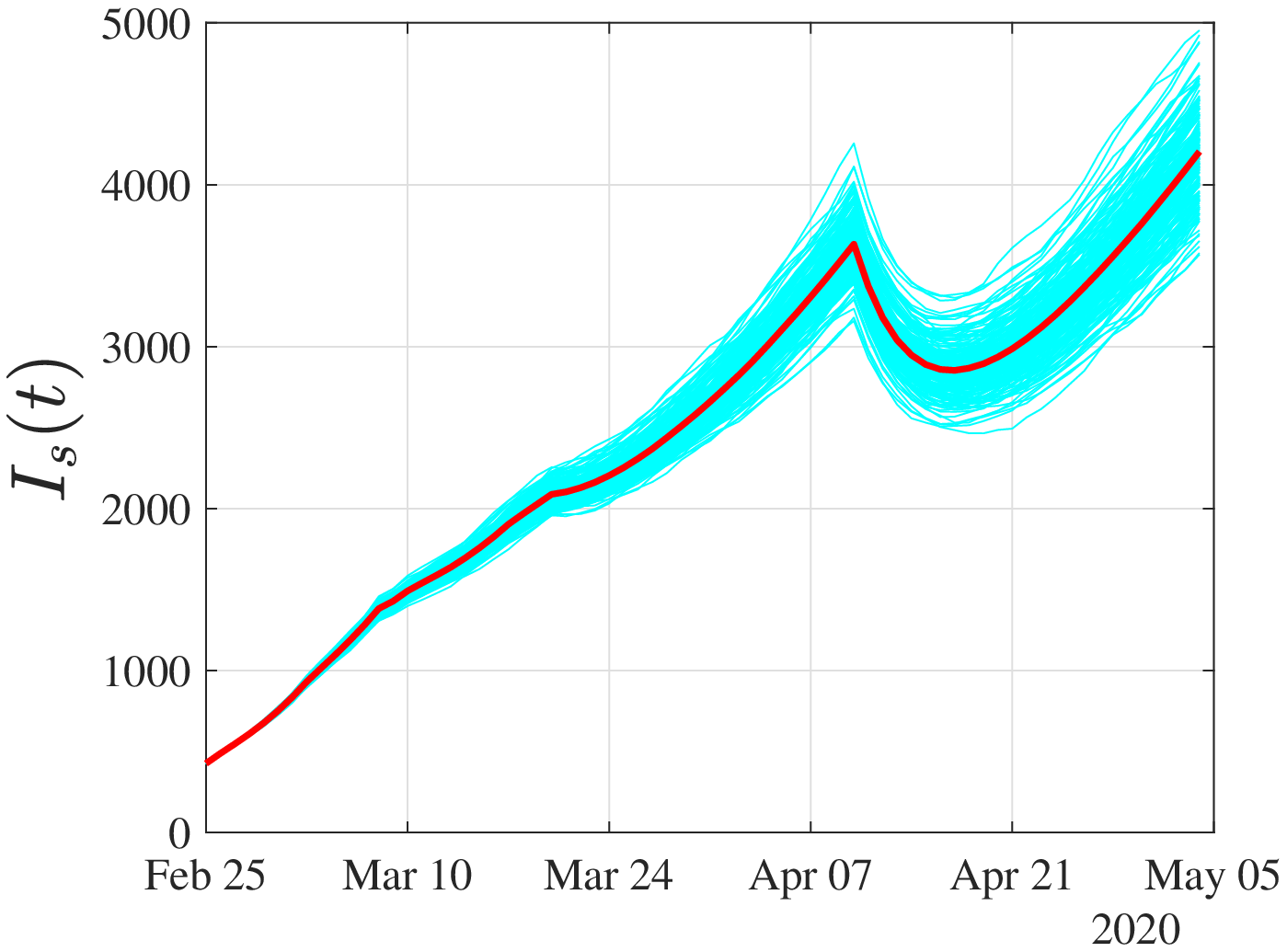}
    }
    \subfigure[]{
         \includegraphics[width=0.31\textheight]{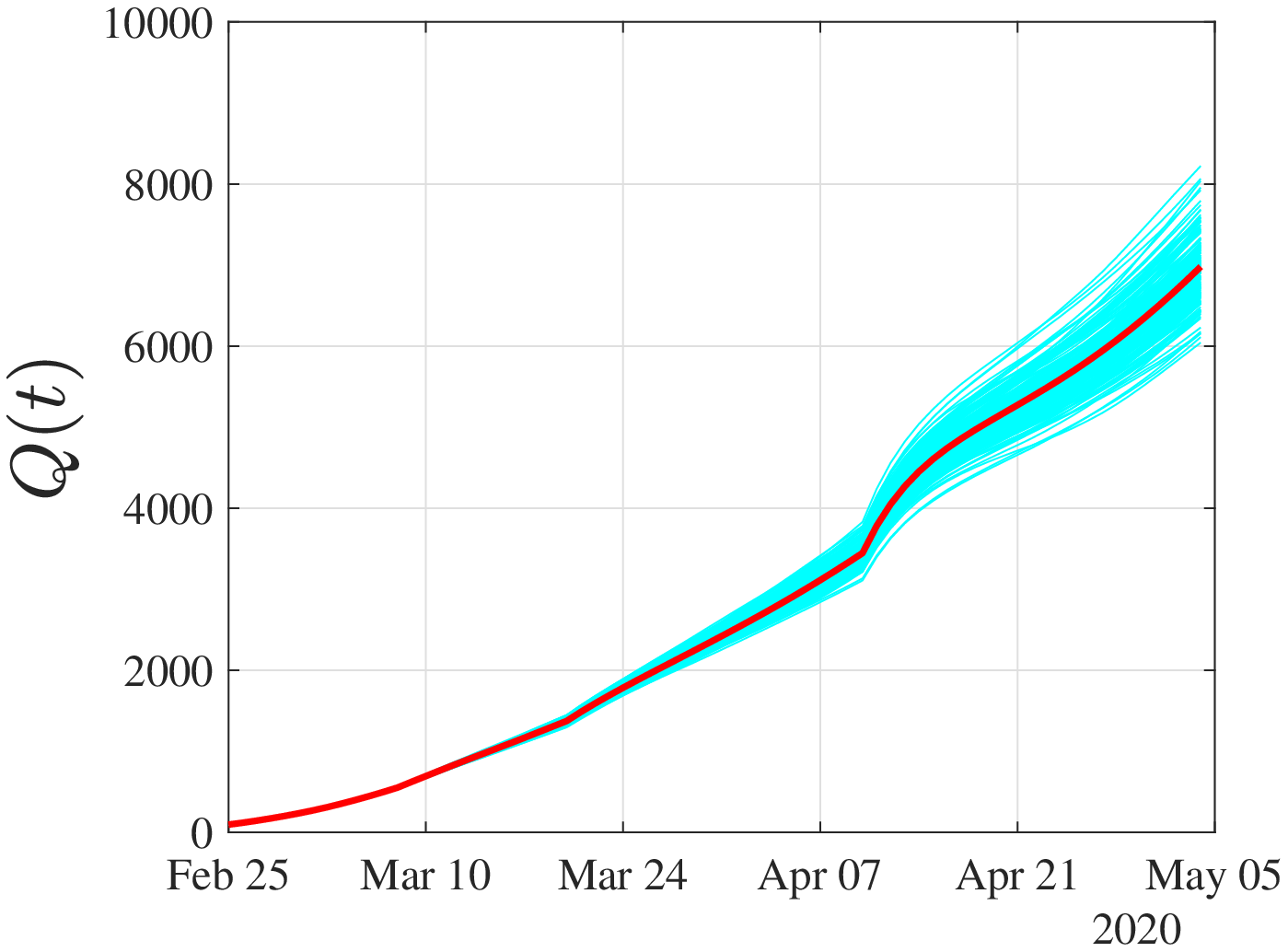}
    }
        \subfigure[]{
         \includegraphics[width=0.31\textheight]{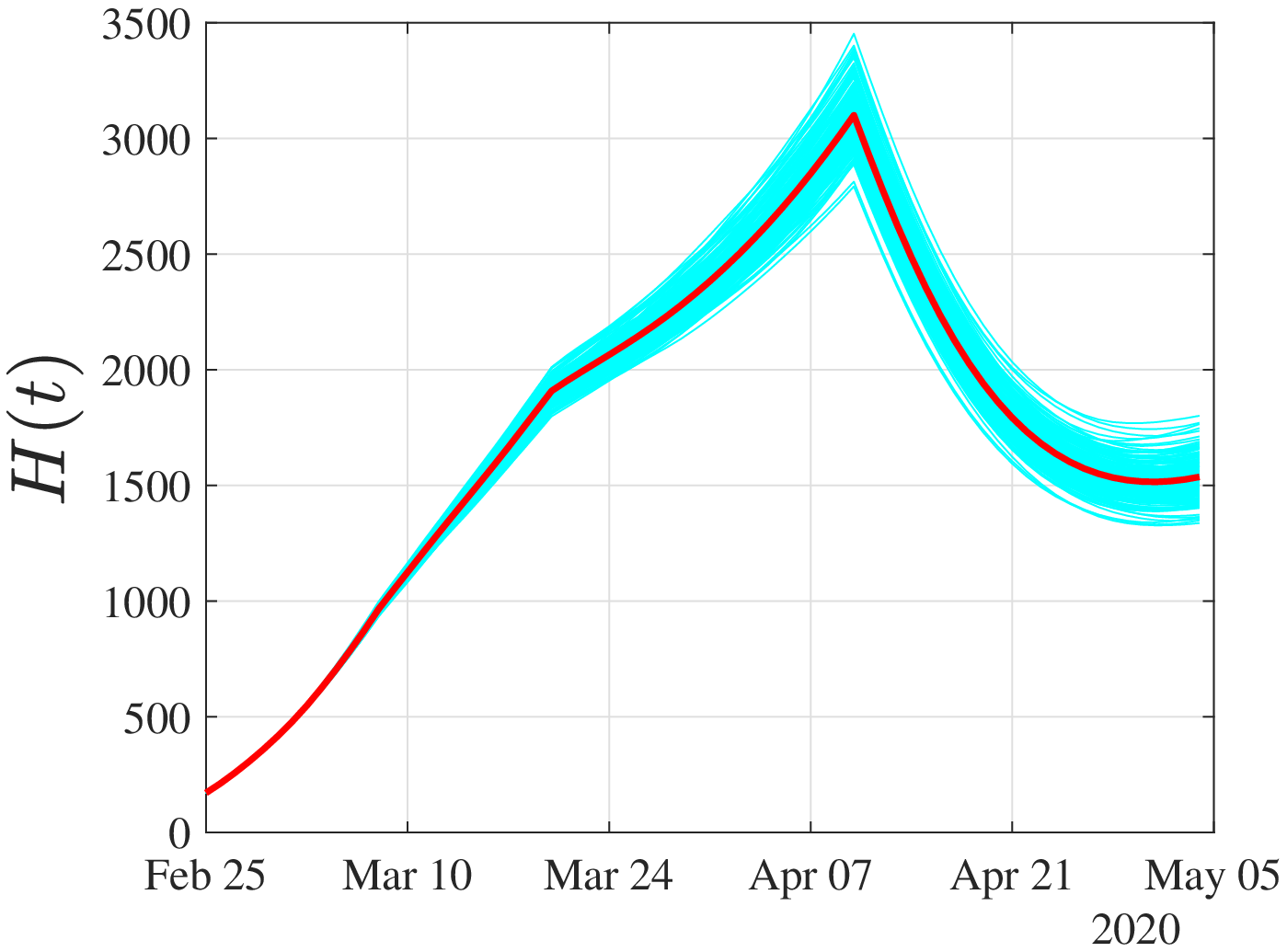}
    }
    \subfigure[]{
         \includegraphics[width=0.31\textheight]{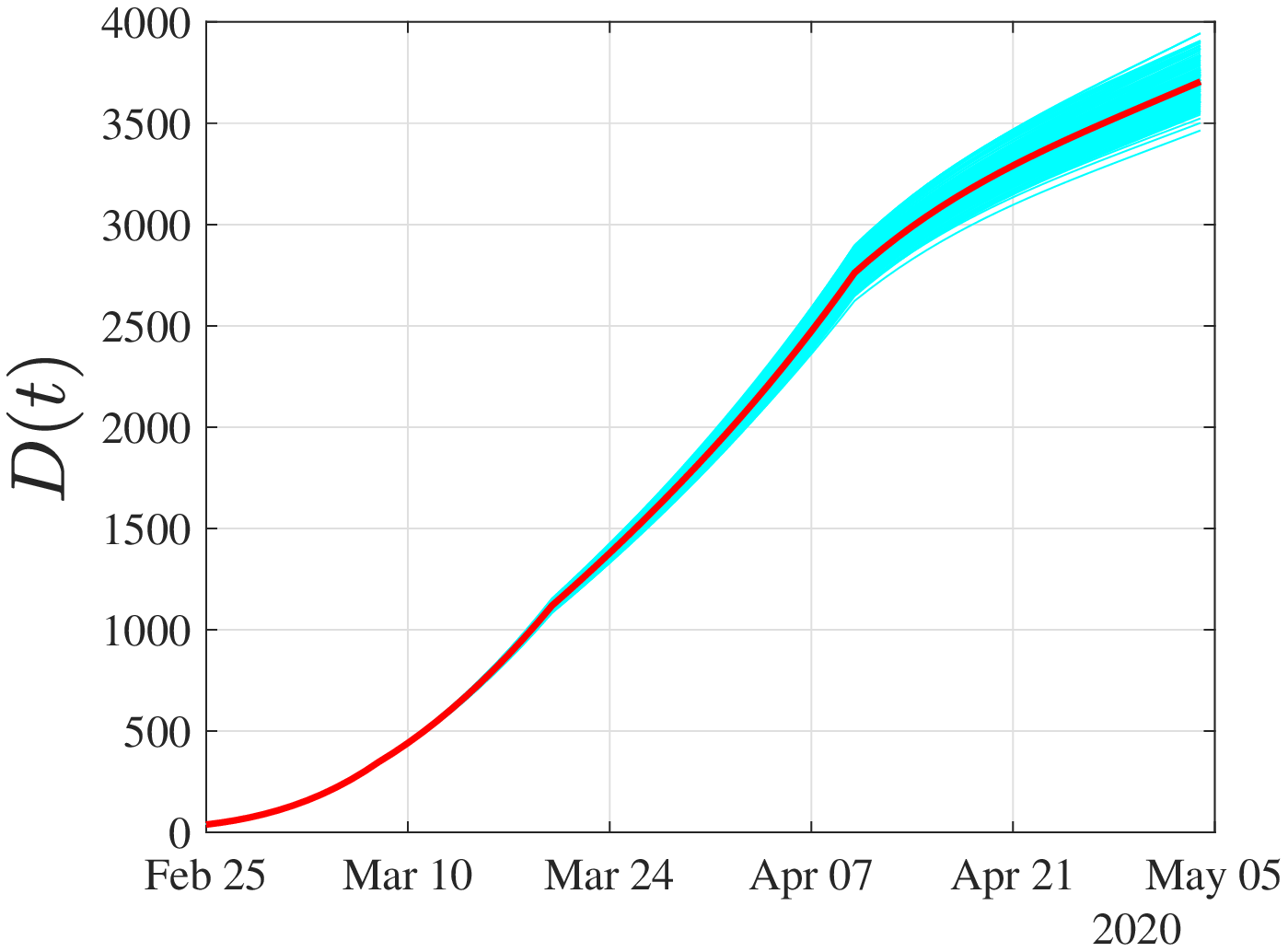}
    }
    \subfigure[]{
         \includegraphics[width=0.31\textheight]{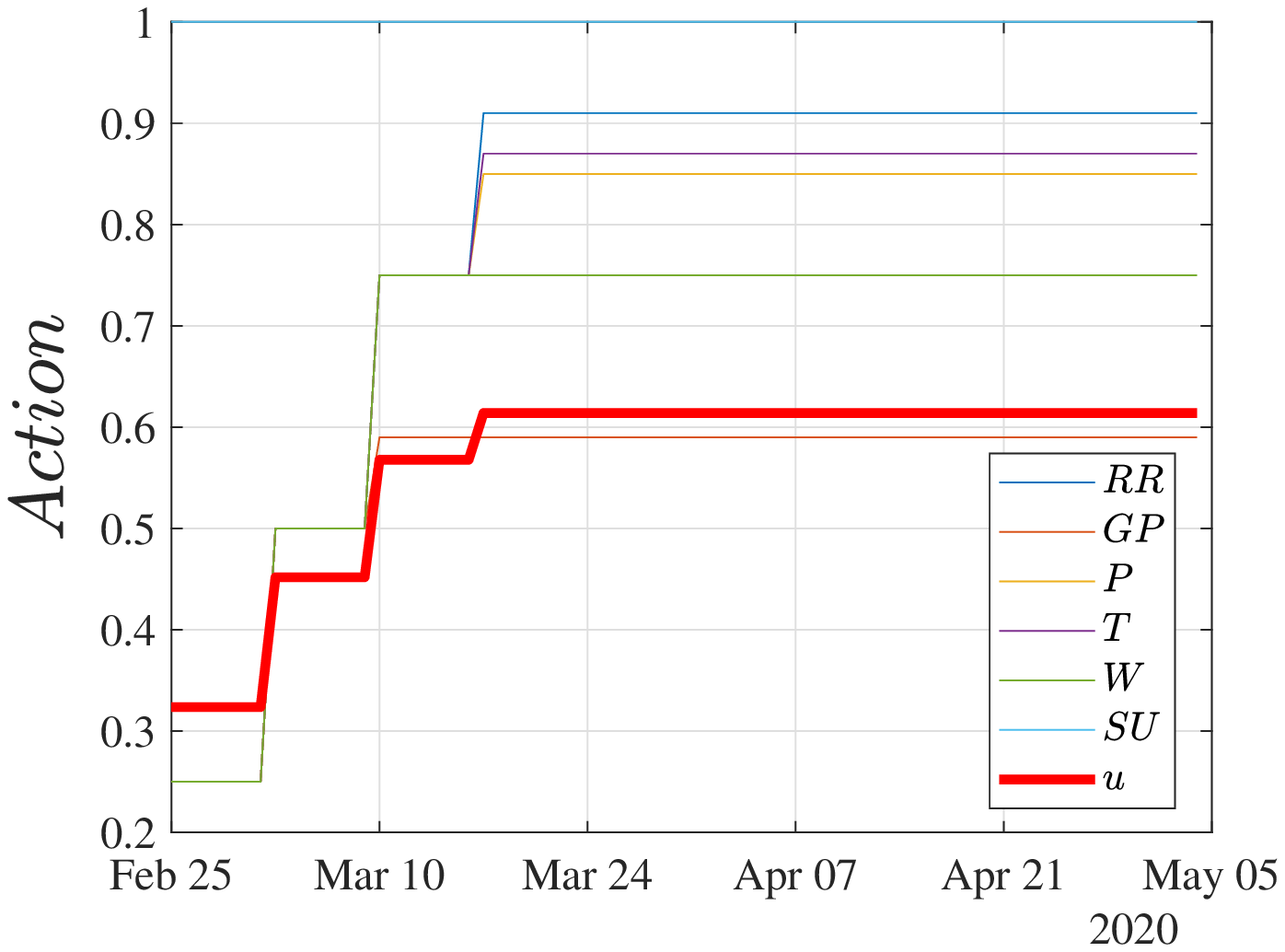}
    }
    \subfigure[]{
         \includegraphics[width=0.31\textheight]{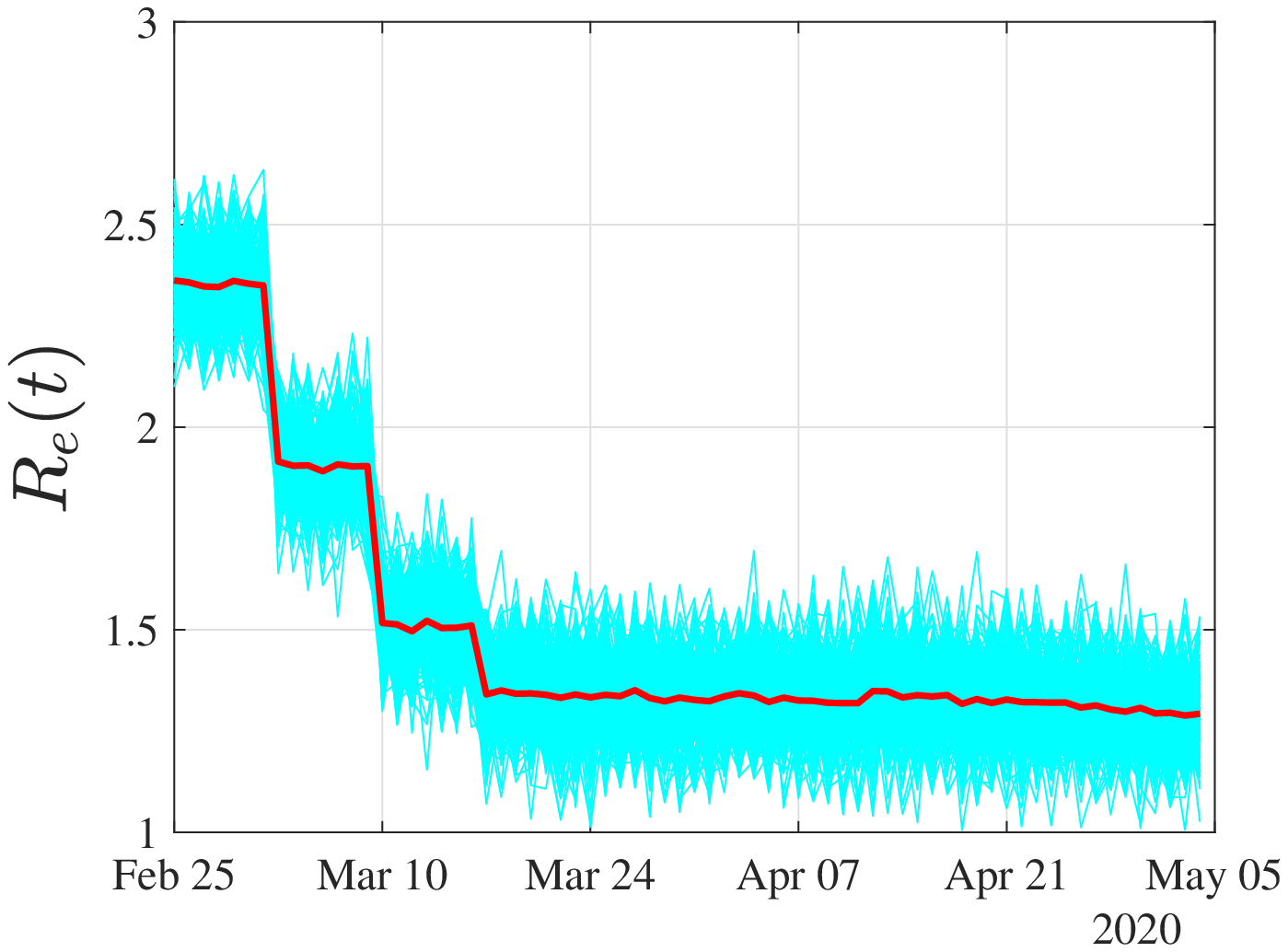}
    }
    \caption{Simulation results with PMPC using the SEASQHRD model for the period February 24-May 4. (a) Infected cases. (b) Hospitalized. (c) Quarantined. (d) Dead. (e) Social distancing in the various activities resulting from the PMPC. (f) Effective Reproduction number $R_e$. Red solid lines in (a), (b), (c), (d) \& (f) correspond to the trajectory expectation. The red solid line in (e) presents the cumulative mobility curtails. Cyan colors present individual trajectories of 200 adherence scenaria.}
    \label{f:MPC_1}
\end{figure}

The choice of the objective function determines the optimization calculations performed by the controller. A linear objective function of \eqref{OP4:J} could be used, but the quadratic objective function of \eqref{OP1:J} resulted in stable control pattern and was used in this work. The probabilistic form of this optimization function was used to account for inherent uncertainty in policy implementation; it is less aggressive than a robust formulation of \eqref{OP3:J}-\eqref{OP3:3}, which assumes the worst-case scenario during optimization. Consequently, the choice of a PMPC formulation allows for less heavy-handed policies while still considering the uncertainty associated with population adherence. 
Another consideration is the relative cost assigned to mobility restrictions, scaled by the relation between parameters $\Omega_A$ and $\omega_I$ \& $\omega_R$. Focusing first on the relative cost of activity restrictions as reflected by the Google mobility report data during the period Feb-24 to May-4, was used to elucidate the preference given in the restrictions. Specifically, the need for groceries is central for a society where activities such as restaurants are curtailed. Similarly, transit is needed in order for basic social activities to take place. On the other hand, parks are not explicitly needed for society to function during epidemic conditions. Using similar arguments a
convenient diagonal matrix $\Omega_A$ was defined, with the diagonal elements expressing the relative resistance to curtails being $\omega_A=[0.2\;1.0\;0.5\;0.5\;1.0\;0.5]$. %In comparison $\Omega_A={\rm\bf I}_{6\times 6}$ 
For a prediction horizon of 3 weeks, control horizon of 2 weeks, and action horizon of 1 week, a selection of $\|\Omega_A\|_F\rightarrow 0$ results in government policies  of aggressive lockdown. A selection of $\|\Omega_A\|_F\approx \omega_R$  results in moderate mobility rstrictions. To equally weight the impact on hospitalization and impact on disease progression we employed $\omega_I \approx 10 \omega_R$. 
Figure \ref{f:MPC_1} shows the simulation results of disease progression for the period Feb-24 to May-4 under the optimal policy as identified by the MPC formulation of \eqref{OP1:J}-\eqref{OP1:3} with $\omega_R=10$, $\omega_I=10$, $\Omega_A$ as discussed previously, $T_p=21$, $T_c=14$, $T_a=7$ days. The uncertainty in the simulation is captured via presenting 200 population adherence scenaria(i.e., $N_{cases}=200$), informed from Table~\ref{t:Google_mobility}.  These parameters ultimately prescribe a government policy that corresponds to approximately a 60\% reduction in mobility. The difference in outcome between the implementation of the prescribed restrictions with the historical data from Lombardy for reference by comparing Figure~\ref{f:MPC_1} with Figure~\ref{f:GM_1}. The simulation of optimal restrictions leads to less infections, hospitalizations, and deaths than the historical data. A major contributor to the difference in deaths is the date of introduction of major curtails which in our simulation starts on Feb-24 which also has the effect that population caution becomes $\theta_C=20\%$ on that day. The historical data in contrast show the effect of casual curtails starting on Feb-21 with major (and aggressive) curtails starting on March 8; population caution becomes $\theta_C=20\%$ on March 8.

\section{Discussion}

One of the more challenging aspects of addressing the COVID-19 crisis is the scale of response, promptness, and fast adaptation to a number of factors ranging from most importantly the strengthening of the public health system both in staff and infrastructures with a particular emphasis on the enhancement of the prevention and primary care facilities network to the design of vaccination policies and from the design of public awareness campaigns to the implementation of appropriate social distancing policies that are needed to combat efficiently this global threat. In an emergent situation in which the inherent uncertainly that characterizes the reported data and generally our knowledge at least at the first stages of an pandemic, mathematical modeling and control theory play an important role towards the simulation and investigation of what if scenarios.

Here, retrospectively, we attempted to assess the impact of social distancing measures on the various activities on the evolution of the pandemic in Lombardy during the first wave of the pandemic (Feb 21-May 4 2020) and its consequence on the fatalities and the pressure on the hospital structures, following a model predictive control approach. Towards this aim, we addressed a compartmental modelling approach that  quantifies the uncertainty of the level of asymptomatic cases \cite{Volpert_2020,Russo2020} in the population and also the effect of social contacts in different activities  to the spread of infectious diseases, In particular, the proposed model extends the one proposed by Russo et al. \cite{Russo2020} for the region of Lombardy, thus considering the hospitalized and guaranteed cases. Furthermore, the effect of the social contacts in various activities is modelled based on the Google mobility reports for the region of Lombardy as well as by metanalyses of the social contacts and their effect in the spread of the infectious diseases as reported for Italy by Mossong et al. \cite{mossong_social_2008}. Thus, the proposed PMPC scheme models also the effect of the reported activities, namely, workplaces, grocery and pharmacy, retail and entertainment, parks and transit to the residential change of mobility. This was achieved by fitting a linear regression model based on the corresponding Google mobility reports.
We show that the proposed modeling framework approximates quite well the actual evolution of the pandemic in terms of reported infected cases, deaths and hospitalized cases. By activating the PMPC scheme from February 21, 2020, i.e. the day that the first official case was reported in Lombardy, we show that the obtained number of deaths would be significantly less (around 4,000 at May 4 instead of the actual number of 14,000), however the spread of the disease would not have been eradicated by May 4. Indeed our simulations show that due to the presence of a significant portion of asymptomatic cases which transfer the disease but are rarely reported, due to the mild or no symptoms, the effective reproduction number is still above 1 on May 4. However, the pressure in the hospitals would be much less, with a pick at a maximum of 3500 hospitalized cases around the mid of April, instead of the actual 13000 cases in the same period. 
Of course, as our study is about a retrospective analysis, it should be regarded as indicative of what would happen if we knew a-priori the scale and severity of the threat. Thus, our analysis should not be considered by any means as a ``ground truth" of any kind of accurate prediction of what would happen if such measures would have been adopted. However,  the results are indicative and highlight the importance of the need of very prompt responses and alertness for analogous emergent and re-emergent infectious diseases that will certainly occur into the future.

On the other hand and importantly, our study pinpoints the urgent need of strengthening the public structures including the epidemiological surveillance mechanisms that will allow faster responses, the enhancement of the (fragmented in most of the countries of the world) primary care facilities networks, the investment in better both qualitatively and quantitatively public hospitals and hiring of medical staff but also in other structures such as schools and mass public transportation. 

The proposed framework can be extended at various fronts. One issue at the forefront of current policy making is vaccination, which this model does not address. To accurately forecast policy decisions today, the model would have to be adjusted for vaccine-induced immunity to COVID-19. Other considerations include the probability of reinfection with the emergence of the virus variants, which was a less significant factor in this study due to the time period examined. Furthermore, the proposed approach could be used for the quantification of the uncertainty of the evolving dynamics, taking into account the reported, from clinical studies, distributions of the epidemiological parameters rather than their expected values and a more detailed compartmental model that takes also the age distribution in the population (see e.g. \cite{blyuss2021effects}. A critical point that is connected with the above is that with such a small number of infectious cases at the initial stage of the simulations, a stochastic model or a hybrid stochastic model could be  more realistic in which uncertainty could be modelled in the form of realistic perturbations (see for example \cite{Onofrio2013}).

Another challenge posed by the COVID-19 pandemic is the accurate state estimation. Here, in order to set the values of the model parameters, we used the reported values of the various epidemiological parameters. However, to examine longer periods of time, methods like moving horizon estimation and design of stochastic state observers, could be used to more accurately capture and assess the uncertainty in the state of the pandemic.
Finally, one could extend the proposed framework by extended as a network of regions and interventions (see \cite{della2020network}. This is an important aspect that we aim to tackle in a future work.

\section*{Acknowledgments}

The work was supported by the Ministero dell'Università e della Ricerca through the program of Fondo Integrativo Speciale per la Ricerca (FISR).Program: FISR2020IP-02893.

\section*{Conflict of interest}

The authors declare there is no conflict of interest.

%\bibliographystyle{unsrt}
%\bibliography{references}  %%% Uncomment this line and comment out the ``thebibliography'' section below to use the external .bib file (using bibtex) .

%\nocite{*}% Show all bib entries - both cited and uncited; comment this line to view only cited bib entries;

\clearpage
\section*{Supporting Information}
%\section{Supporting Information}
\setcounter{section}{0}
\renewcommand*{\thesection}{S.\arabic{section}}
\setcounter{figure}{0}
\renewcommand{\thefigure}{S\arabic{figure}}
\setcounter{table}{0}
\renewcommand{\thetable}{S\arabic{table}}

\section{Epidemiological Data}

All the relevant data used in this paper are publicly available and accessible at \url{https://lab.gedidigital.it/gedi-visual/2020/coronavirus-i-contagi-in-italia/}. The reported cumulative numbers of cases from February 21 to March 19 are listed in Table S1. The data from February 21 to March 8 have been used for the calibration of the model parameters and the data from March 9 to March 19 have been used for the validation of the model.

\section{Google mobility Data}
All the mobility data used in the present study are publicly available at https://www.google.com/covid19/mobility/.
These data were employed to estimate the standard deviation and bias between the government guidelines and population behavior, discussed in \S\ref{s:Mobility}.

%%%%%
\section{Other MPC formulations}\label{s:MPCs}
Beyond the PMPC formulation of \eqref{OP1:J}-\eqref{OP1:3}, deterministic and robust MPC problems, (\eqref{OP2:J}-\eqref{OP2:3} and \eqref{OP3:J}-\eqref{OP3:3} respectively, were also formulated and their performance was evaluated.

The deterministic problem does not consider population adherence to government guidelines, assuming perfect response. This leads to a simpler formulation at the cost of a potentially optimistic view of the progression of the epidemic that can lead to the violation of constraints \eqref{OP2:3} and importantly of \eqref{OP2:2}.
\begin{align}
&    \alpha^*(t)=\arg\min_{\delta \alpha\in\mathbb{R}^{6\times N_d}}
\sum_{t=\text{Day 1}}^\text{Day $T_P$} \omega_I\left(\frac{{\rm \bf E}[H(t;u,\theta_A)]}{{\rm \bf E}[H(t;u\equiv 0,\theta_A)]}\right)^2 
+ \omega_R ({\rm \bf E}\left[ R_{eff}(t;u,\theta_A)\right])^2
+ \alpha(t)^T\Omega_A\alpha(t)\label{OP2:J}\\
&    \text{subject to}\nonumber\\
&    \alpha(t)=\sum_{j=1}^{N_d} {\bf H}(t-(j-1)\,T_A)\delta \alpha_{j}\label{OP2:u}\\
&    {\bf 0}\le \alpha(t)\le [0.9\; 0.6\; 0.8\; 0.9\; 0.8\; 1]^T \label{OP2:0}\\
&    \delta\alpha_j\le \delta\alpha_c \label{OP2:0a}\\
&    u(t)=[0.216\;0.076\;0.04\;0.063\;0.117\;0.196]\alpha(t) \label{OP2:_1}\\
&    \beta(t)=\beta_0(1-\theta_C(t))(1-u(t)+\theta_A(t))\\
&    {\rm \bf E}[H(t; u,\theta_A)] \le Beds, \quad t=\{1,...,T_p\}\label{OP2:1}\\
&    {\rm \bf E}[R_{eff}(t; u,\theta_A)] \le R_{e,c}, \quad t=\{1,...,T_p\}\label{OP2:2}\\
&    X(t;u,\theta_A)= {\rm SEASQHRD}(X(t-1;u,\theta_A),u(t),\theta_A),\;\; X(0; u,\theta_A)=X_0,\; \; t=\{1,...,T_p\}\label{OP2:3}
%&    u(t)=0.216\,RR(t)+0.076\,G(t)+0.04\,P(t)+0.063\,T(t)+0.117\,W(t)+0.196\,S(t) \label{OP1:_1}\\
%    
\end{align}
At the other end of the spectrum, the robust formulation considers  bounded uncertainties in the state evolution and the enforcement of the control action. It attempts to identify a policy that will not violate the constraints for all investigated scenaria. It is important to note that basic robust MPCs are designed to address problems where the uncertainty is bounded; elaborate RMPC formulations have been designed using Lyapunov arguments \cite{1532396,MHASKAR2005209}. In the presented work, this is addressed by creating a set of ``scenaria"; each scenario contains values of the adherence of the population to each activity curtail and the associated expected initial state $X_0$ that will be used by the SEASQHRD model.
A robust optimization problem was also formulated where a non-convex set was defined containing value ranges for each activity and unmeasurable initial state. As the results from these two robust MPC formulations were almost identical, we only present the first RMPC formulation here.
\begin{align}
&    \alpha^*(t)=\arg\min_{\delta \alpha\in\mathbb{R}^{6\times N_d}}\max_{\theta_A\in\Theta}
\sum_{t=\text{Day 1}}^\text{Day $T_P$} \omega_I\left(\frac{H(t;u,\theta_A)}{H(t;u\equiv 0,\theta_A)}\right)^2
+ \omega_R R_{eff}(t;u,\theta_A)^2
+ \alpha(t)^T\Omega_A\alpha(t)\label{OP3:J}\\
&    \text{subject to}\nonumber\\
&    \alpha(t)=\sum_{j=1}^{N_d} {\bf H}(t-(j-1)\,T_A)\delta \alpha_{j}\label{OP3:u}\\
&    {\bf 0}\le \alpha(t)\le [0.9\; 0.6\; 0.8\; 0.9\; 0.8\; 1]^T \label{OP3:0}\\
&    \delta\alpha_j\le \delta\alpha_c \label{OP3:0a}\\
&    u(t)=[0.216\;0.076\;0.04\;0.063\;0.117\;0.196]\alpha(t) \label{OP3:_1}\\
&    \beta(t)=\beta_0(1-\theta_C(t))(1-u(t)+\theta_A(t))\\
&     \max_{\theta_A\in\Theta}H(t; u,\theta_A) \le Beds, \quad t=\{1,...,T_p\}\label{OP3:1}\\
&    \max_{\theta_A\in\Theta}R_{eff}(t; u,\theta_A) \le R_{e,c}, \quad t=\{1,...,T_p\}\label{OP3:2}\\
&    X(t;u,\theta_A)= {\rm SEASQHRD}(X(t-1;u,\theta_A),u(t),\theta_A(t)),\;\; X(0; u,\theta_A)=X_0,\; \; t=\{1,...,T_p\}\label{OP3:3}
%&    u(t)=0.216\,RR(t)+0.076\,G(t)+0.04\,P(t)+0.063\,T(t)+0.117\,W(t)+0.196\,S(t) \label{OP1:_1}\\
%    
\end{align}
Finally, recognizing  that $H$ and $\alpha$ always attain positive values, a linear objective function can be considered, compared to the widely used quadratic one:
\begin{equation}
\displaystyle J_P=
\sum_{t=\text{Day 1}}^\text{Day $T_P$} \omega_I {\rm \bf E}\left[\frac{H(t;u,\theta_A)}{H(t;u\equiv 0,\theta_A)}\right] 
+ \omega_R {\rm \bf E}\left[ R_{eff}(t;u,\theta_A)\right]
+ {\bf 1}_{1\times 6}\Omega_A\alpha(t))\label{OP4:J}
\end{equation}
\begin{equation}
\displaystyle J_E=
\sum_{t=\text{Day 1}}^\text{Day $T_P$} \omega_I\frac{{\rm \bf E}[H(t;u,\theta_A)]}{{\rm \bf E}[H(t;u\equiv 0,\theta_A)]} 
+ \omega_R {\rm \bf E}\left[ R_{eff}(t;u,\theta_A)\right]
+ {\bf 1}_{1\times 6}\Omega_A\alpha(t)\label{OP5:J}
\end{equation}
%%%%%

\section{Effect of MPC parameters on policy design}
\begin{figure}[p]
    \centering
    \subfigure[]{
         \includegraphics[width=0.31\textheight]{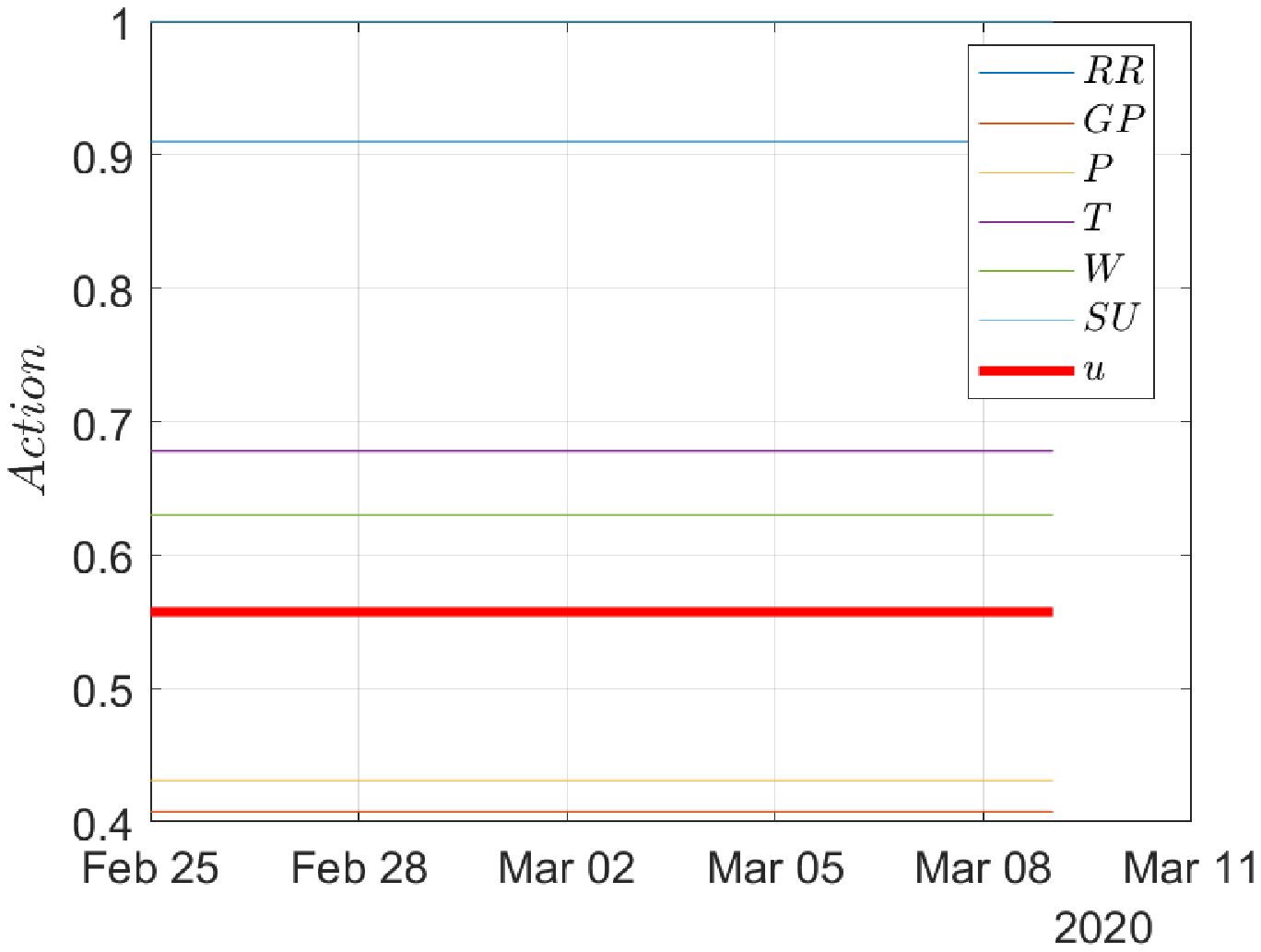}
    }
    \subfigure[]{
         \includegraphics[width=0.31\textheight]{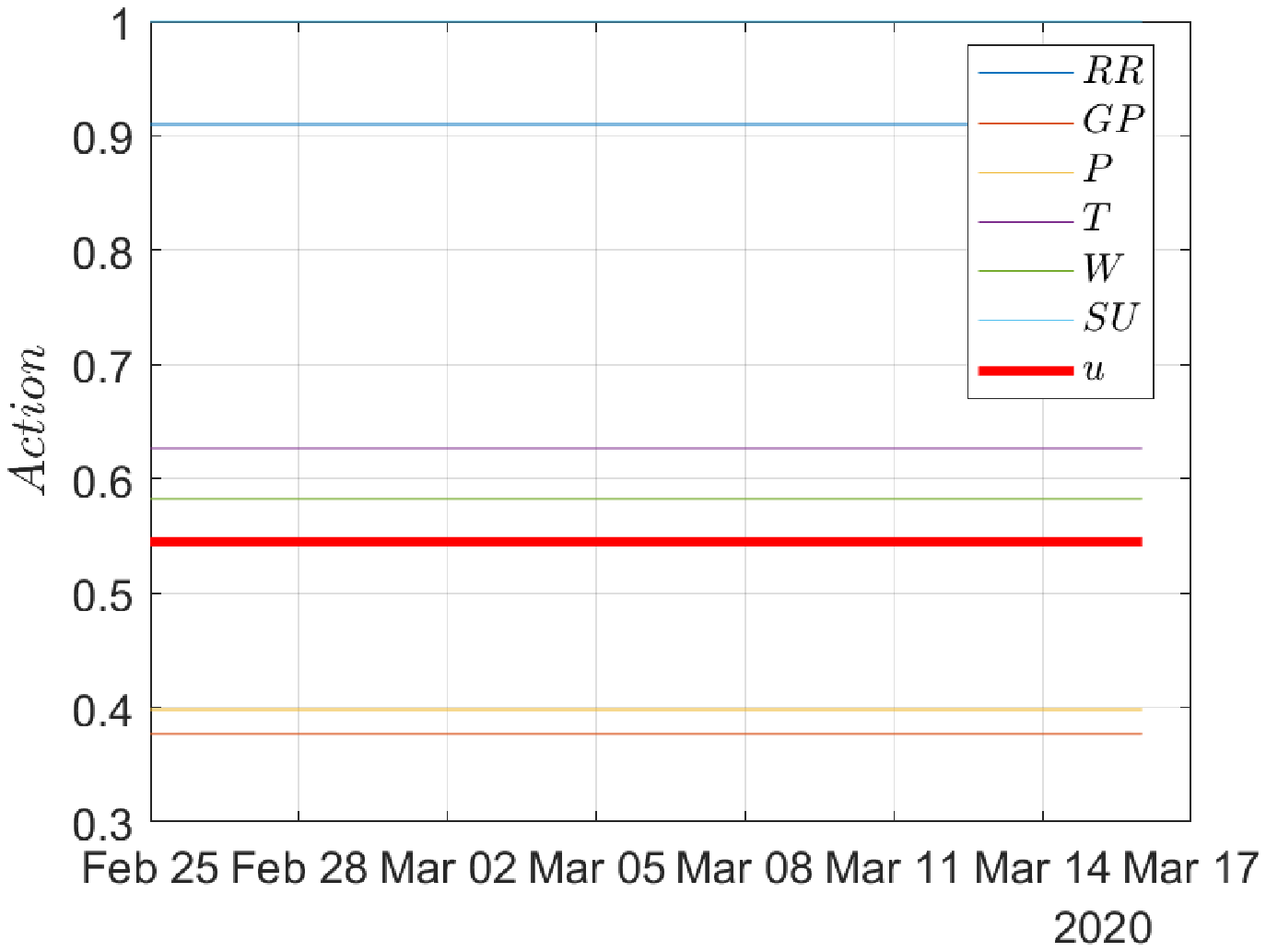}
    }\\
        \subfigure[]{
         \includegraphics[width=0.31\textheight]{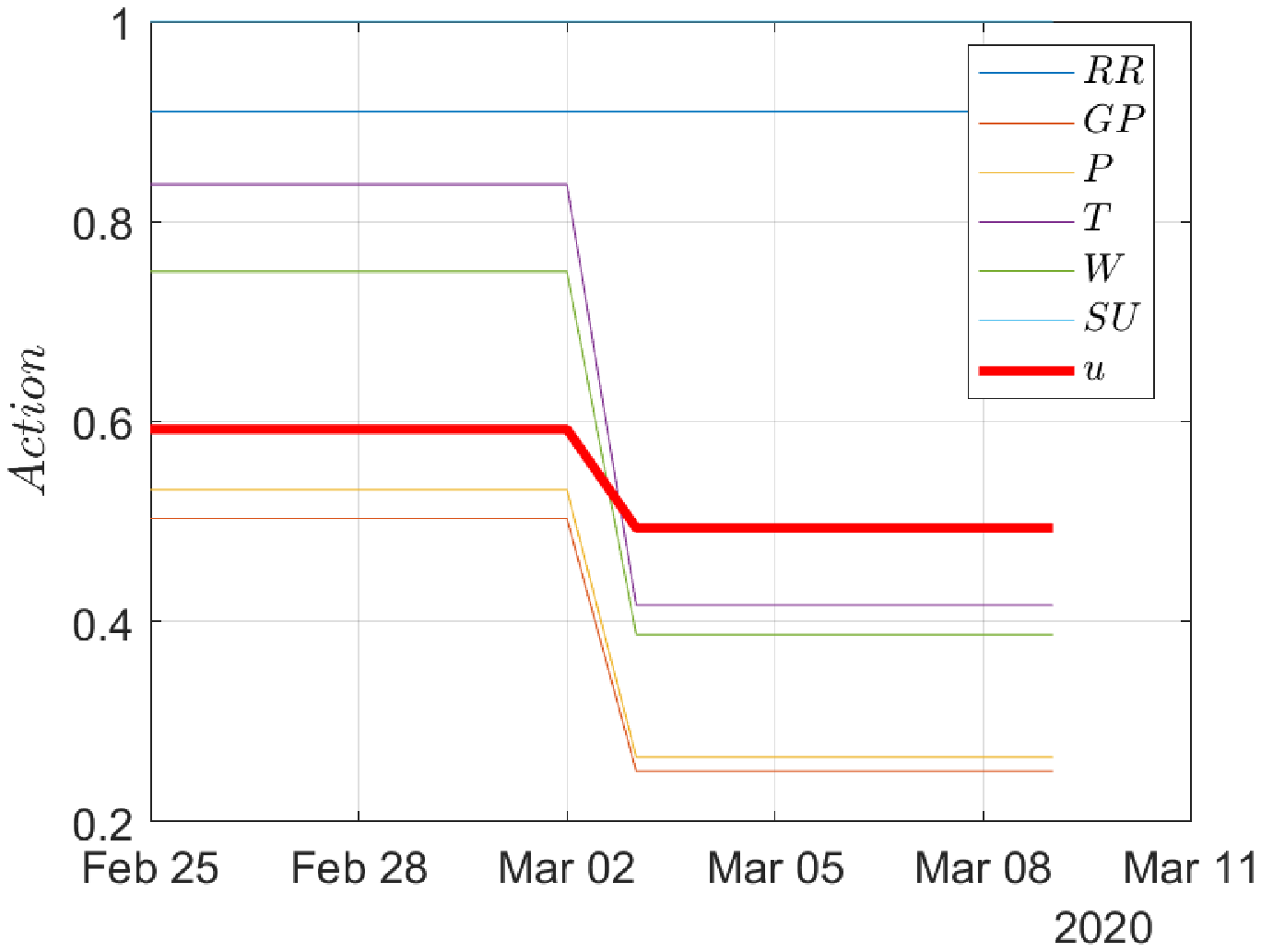}
    }
    \subfigure[]{
         \includegraphics[width=0.31\textheight]{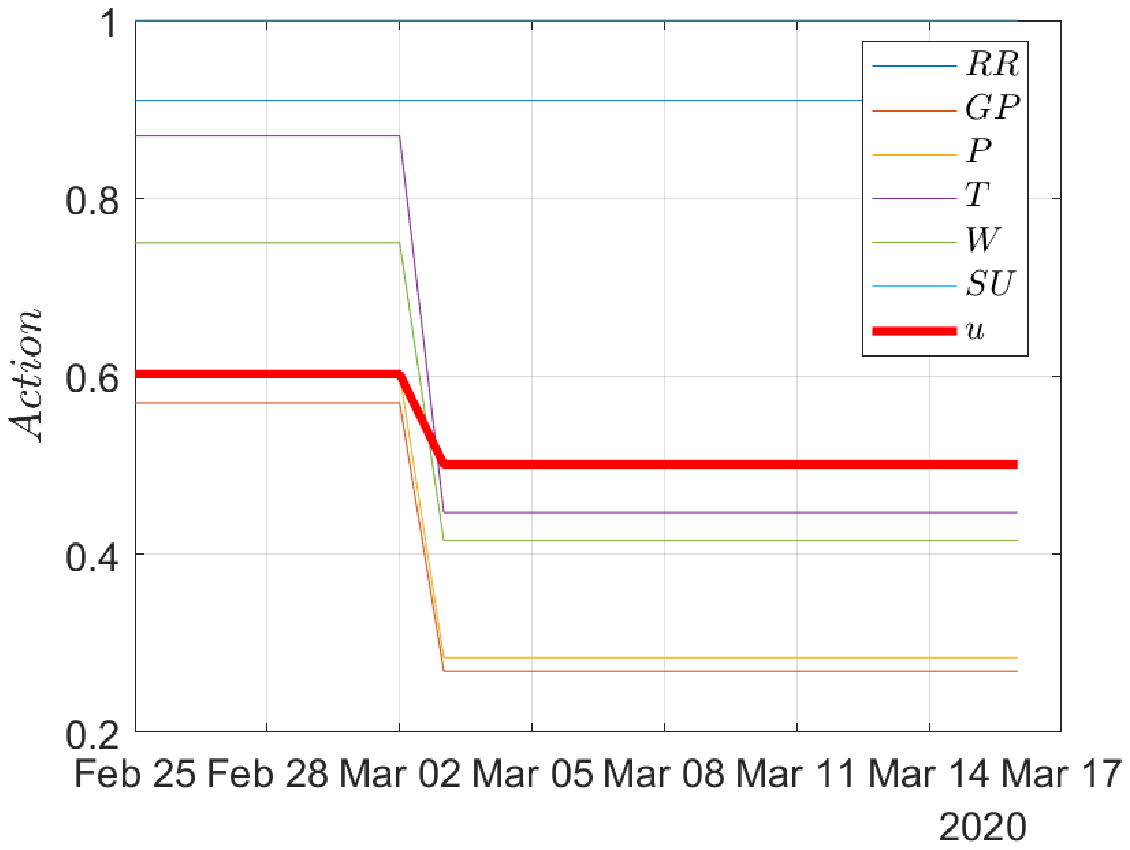}
    }
    \subfigure[]{
         \includegraphics[width=0.31\textheight]{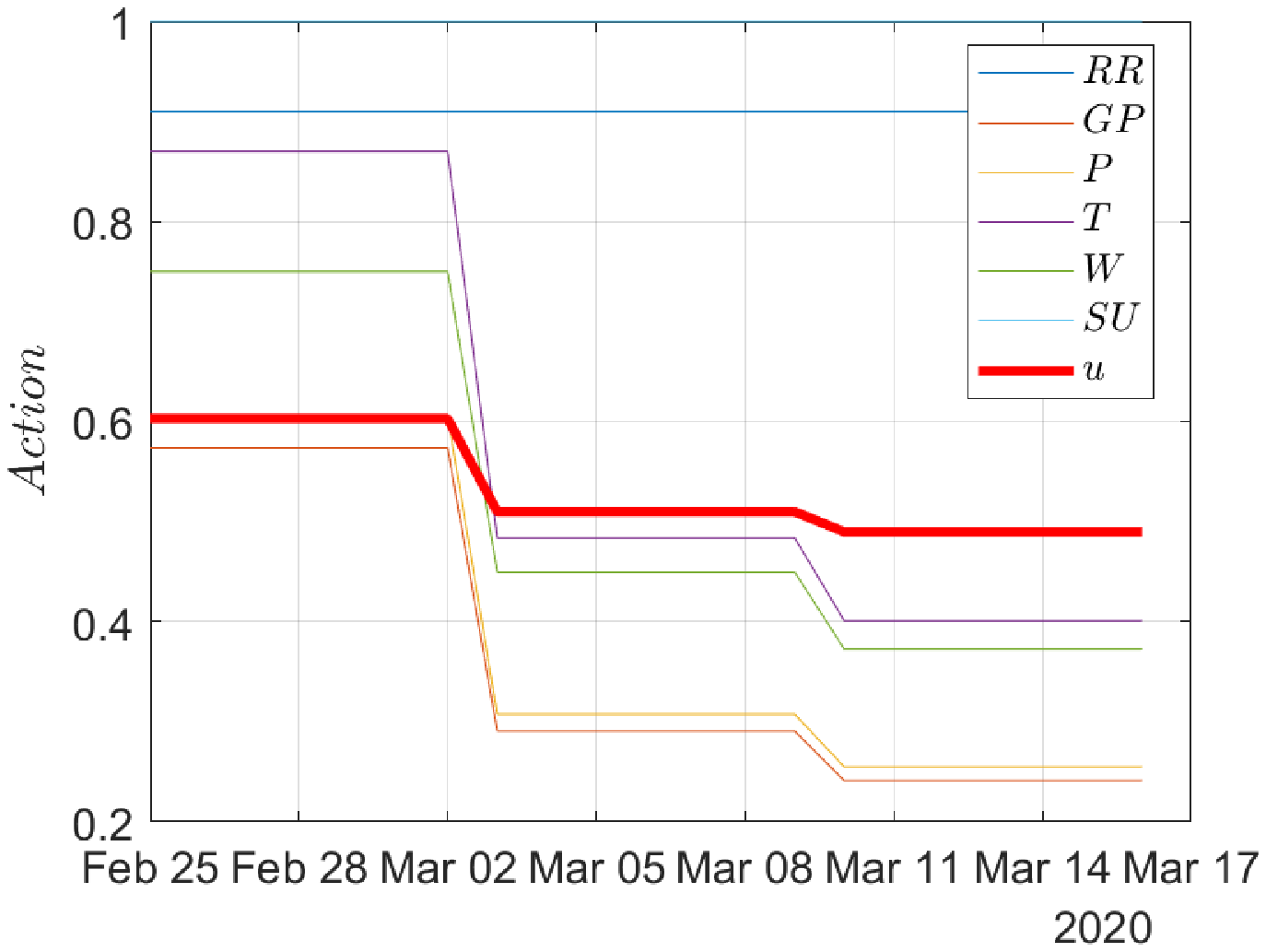}
    }
    \caption{Optimal action from MPC of formulation \eqref{OP1:J}-\eqref{OP1:3} with $\omega_I=10$, $\omega_R=0.5$ and $\Omega_A$ presented in the main manuscript for different horizons: (a) $T_p=14$, $T_c=7$; (b) $T_p=14$, $T_c=14$; (c) $T_p=21$, $T_c=7$; (d) $T_p=21$, $T_c=14$; (e) $T_p=21$, $T_c=21$. Red solid lines  correspond to the aggregate mobility restriction. We observe that as the prediction horizon increases the action of the 1$^{\rm st}$ period converges to a set value.}
    \label{f:MPC_S1}
\end{figure}

Numerous simulations were performed to investigate the effect MPC parameters have on proposed action. Figure ~\ref{f:MPC_S1} showcases the effect of the MPC horizon parameters on the calculated optimal action. We chose horizons $T_p=21$, $T_c=14$, $T_a=7$ days that balanced problem complexity and accuracy well.

\begin{figure}[p]
    \centering
    \subfigure[]{
         \includegraphics[width=0.31\textheight]{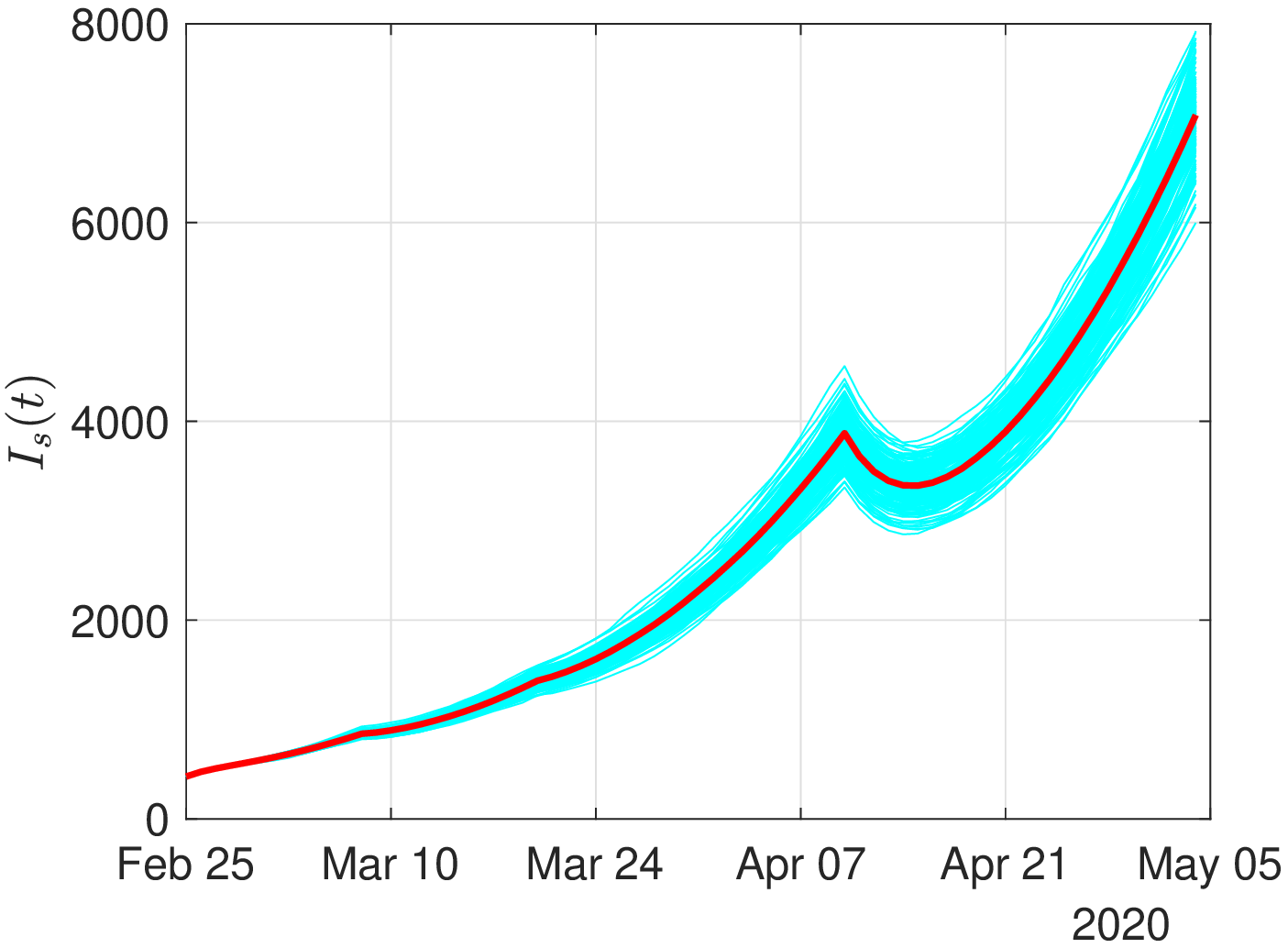}
    }
    \subfigure[]{
         \includegraphics[width=0.31\textheight]{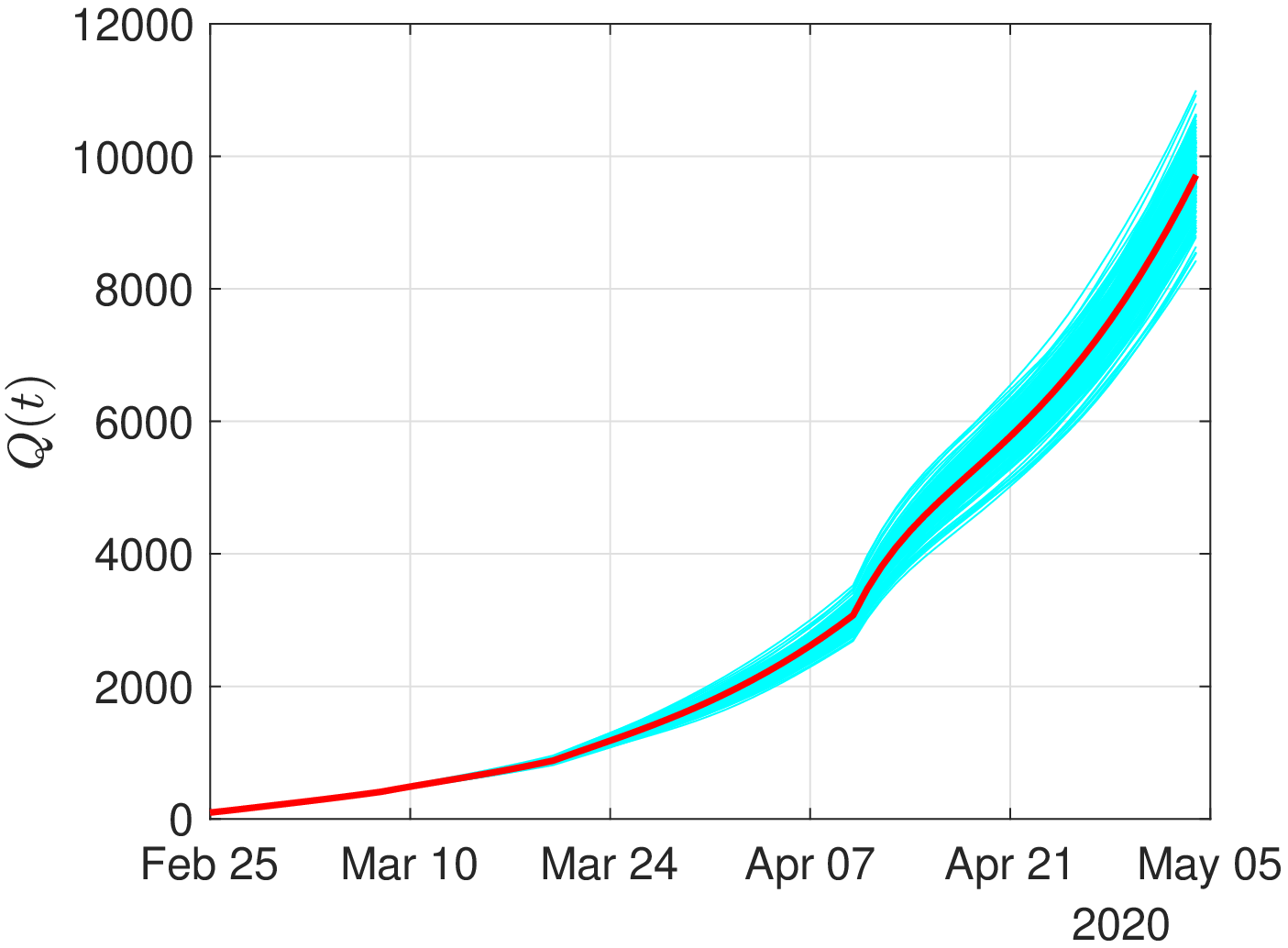}
    }
        \subfigure[]{
         \includegraphics[width=0.31\textheight]{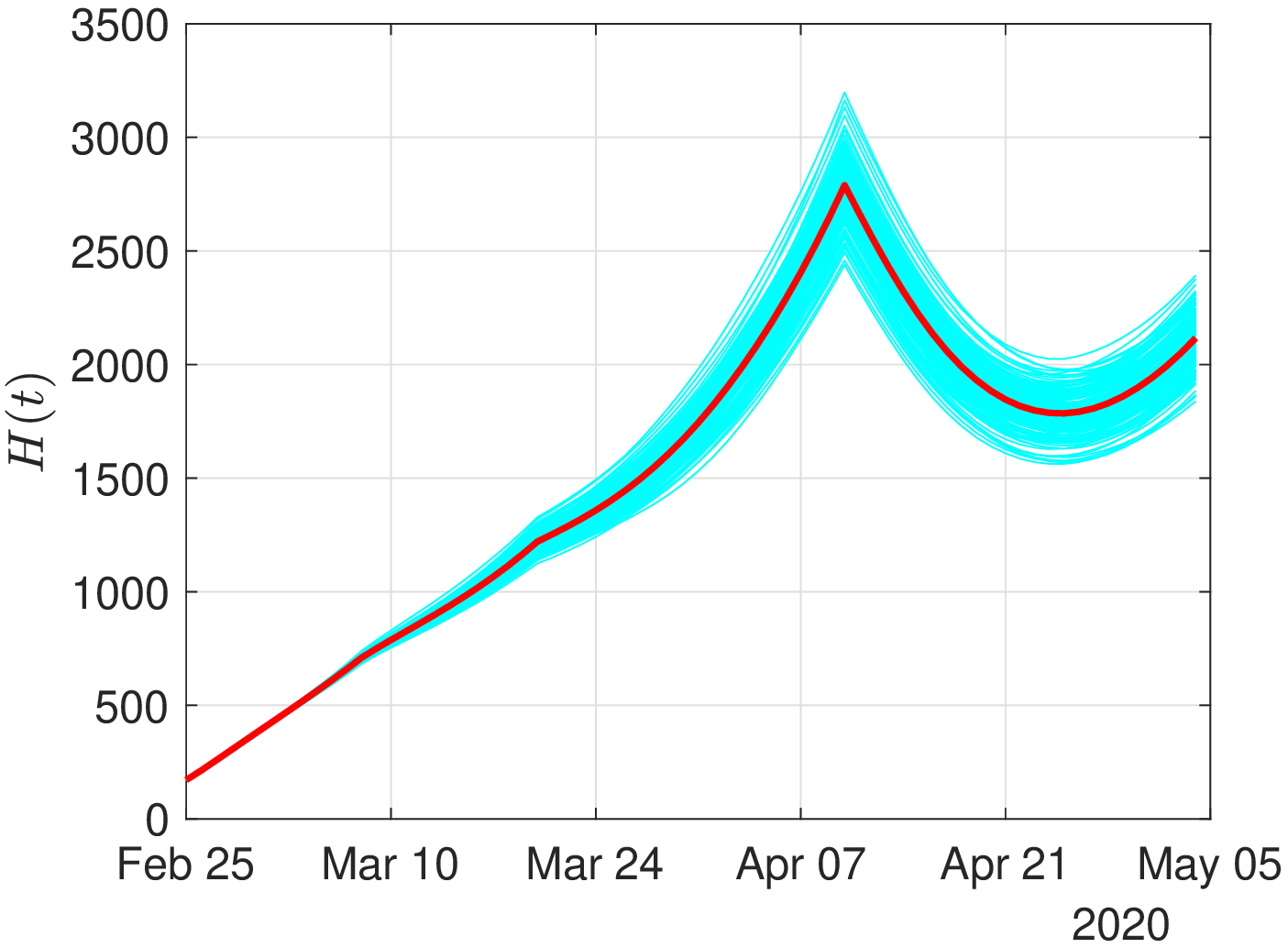}
    }
    \subfigure[]{
         \includegraphics[width=0.31\textheight]{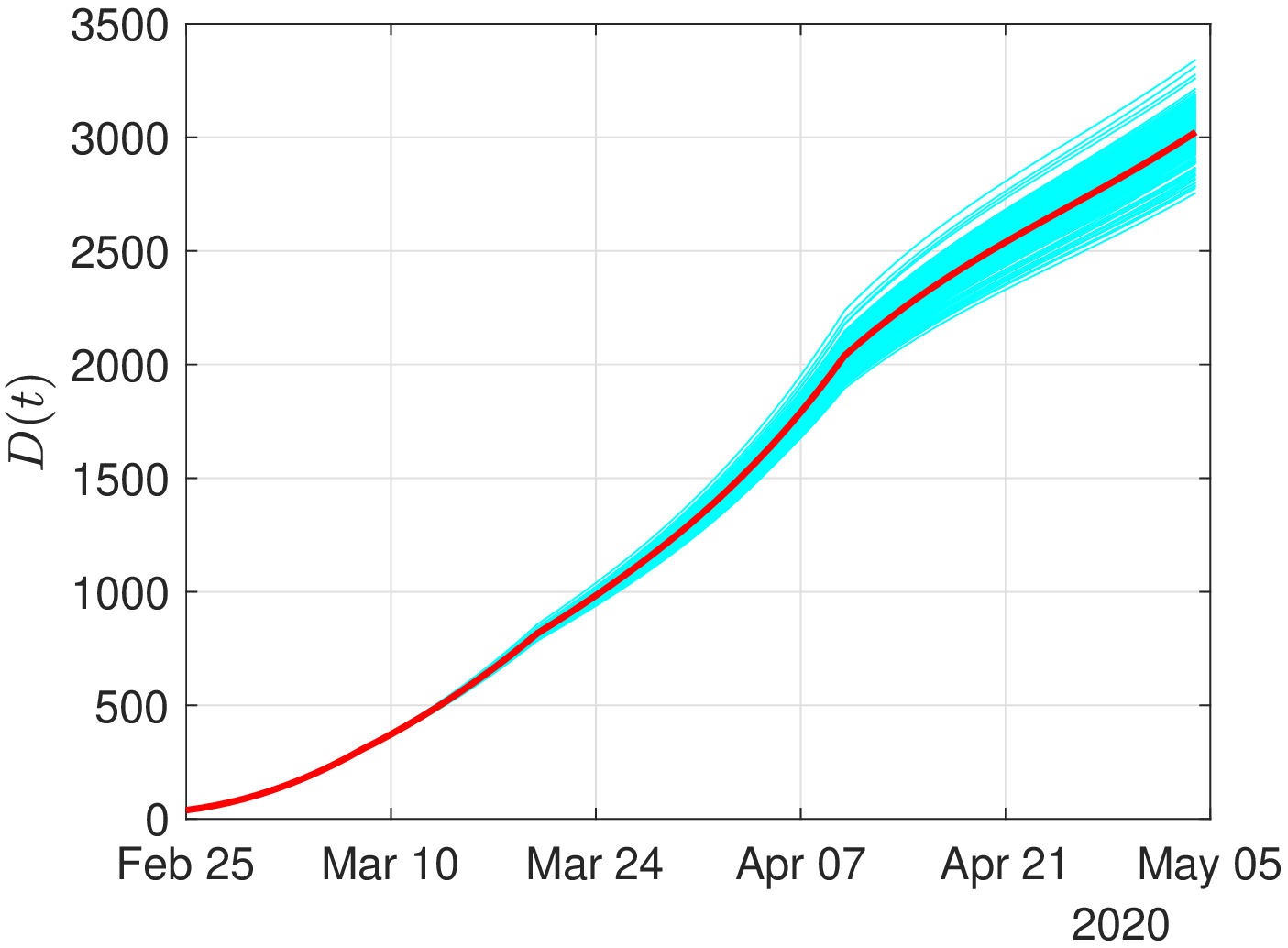}
    }
    \subfigure[]{
         \includegraphics[width=0.31\textheight]{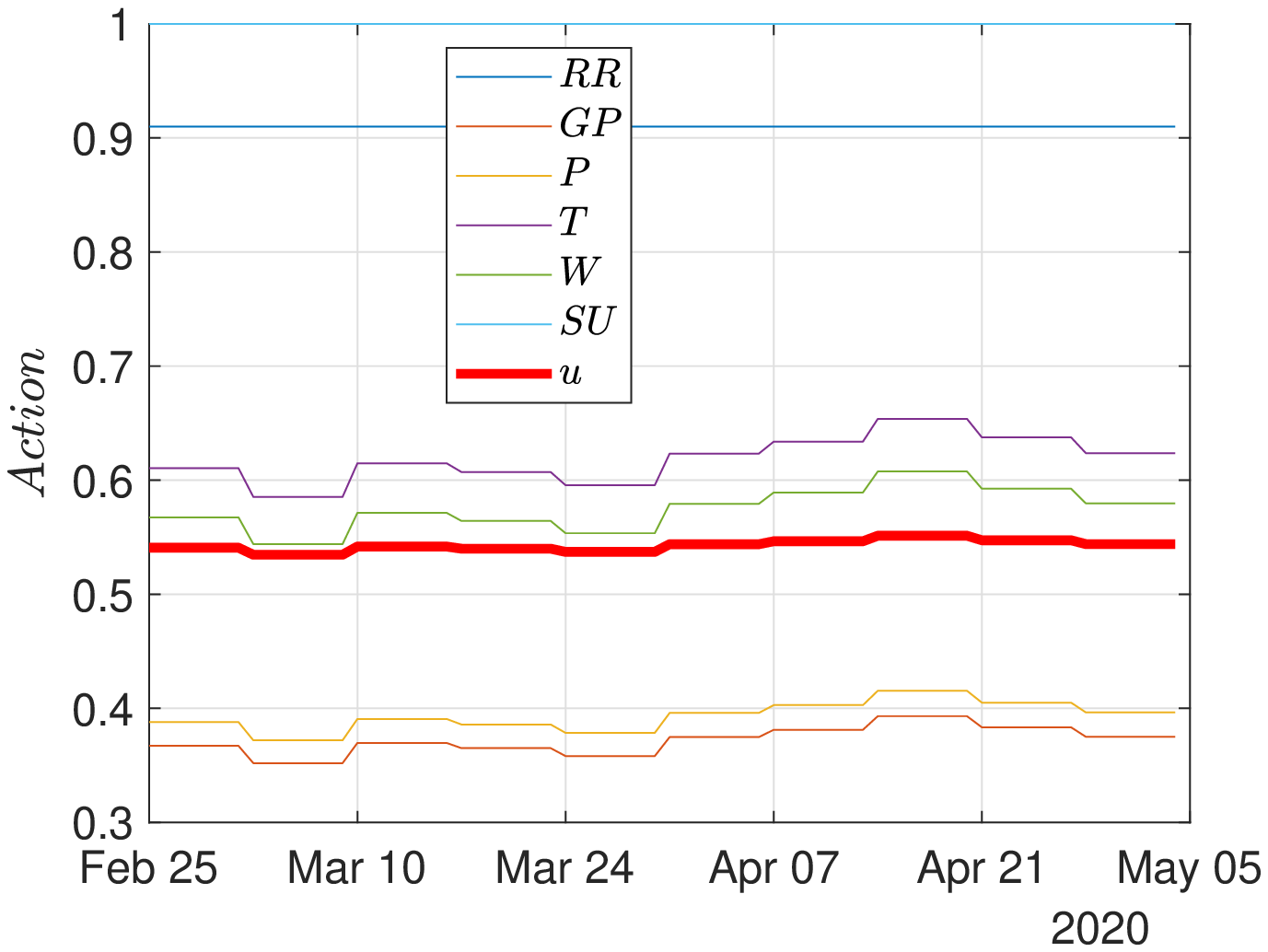}
    }
    \subfigure[]{
         \includegraphics[width=0.31\textheight]{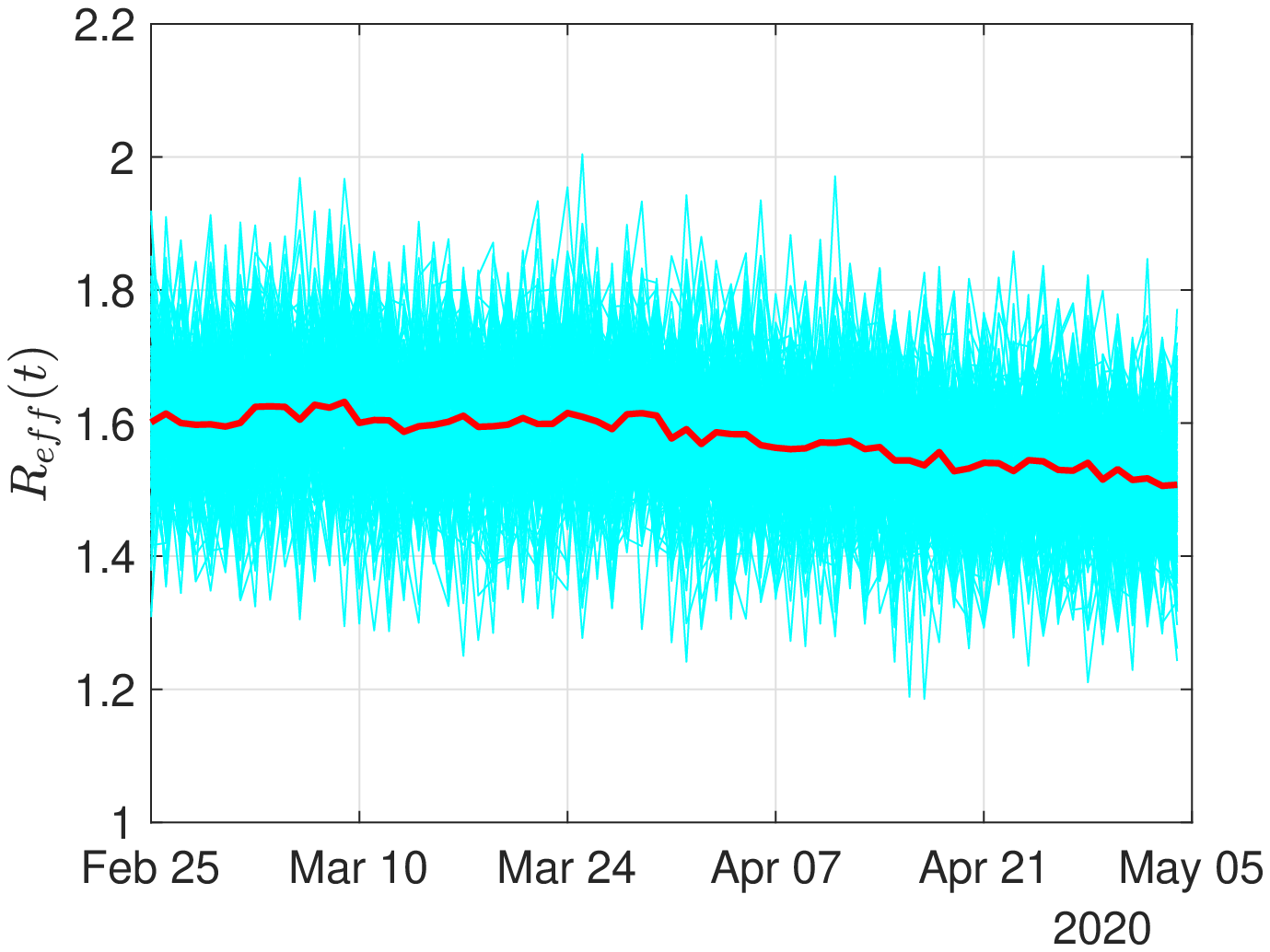}
    }
    \caption{Simulation results of epidemic progression under MPC for the period February 24-May 4 with $\omega_I=0$, $\omega_R=1$, $\Omega_A$ as in \eqref{e:V_A} and $\delta\alpha_c={\bf 1}_{6\times 1}$ in \eqref{OP1:0a}. (a) Infected cases. (b) Hospitalized. (c) Quarantined. (d) Dead. (e) Social distancing in the various activities resulting from the MPC. (f) Effective Reproduction number $R_e$. Red solid lines in (a), (b), (c), (d) \& (f) correspond to the trajectory expectation. The red solid line in (e) presents the cumulative mobility curtails. Cyan colors present individual trajectories of 200 adherence scenaria.}
    \label{f:MPC_S2}
\end{figure}

\begin{figure}[p]
    \centering
    \subfigure[]{
         \includegraphics[width=0.31\textheight]{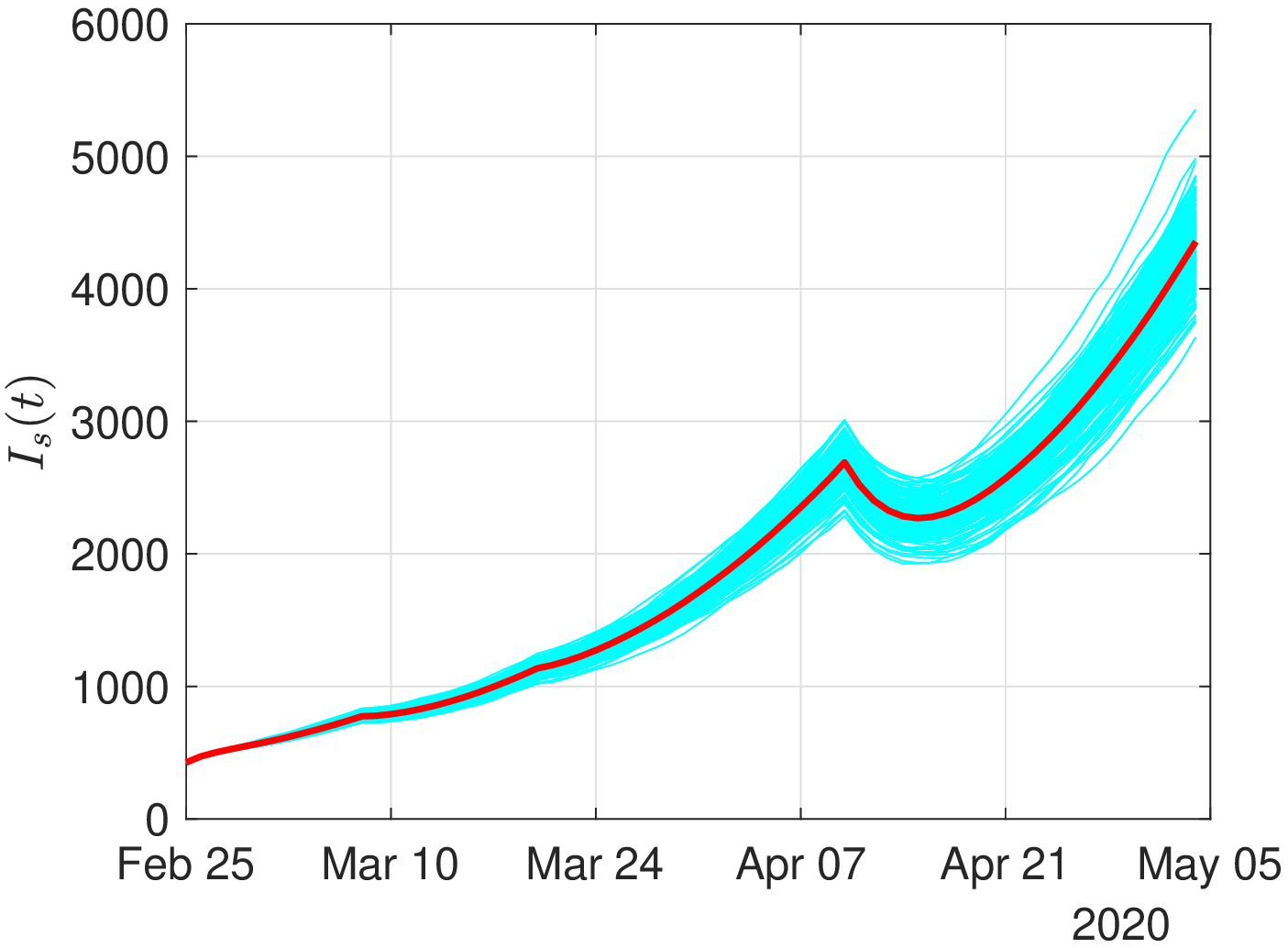}
    }
    \subfigure[]{
         \includegraphics[width=0.31\textheight]{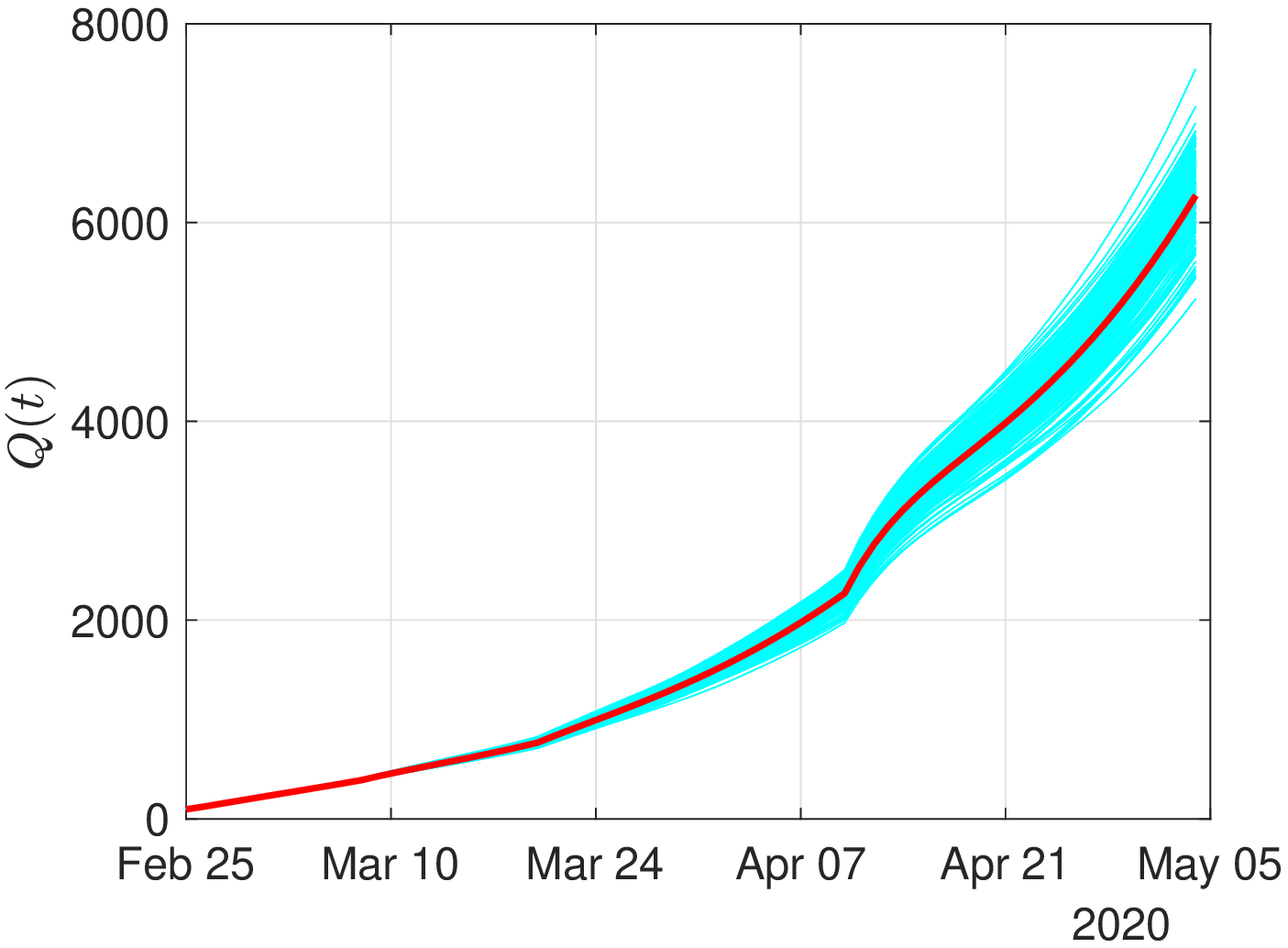}
    }
        \subfigure[]{
         \includegraphics[width=0.31\textheight]{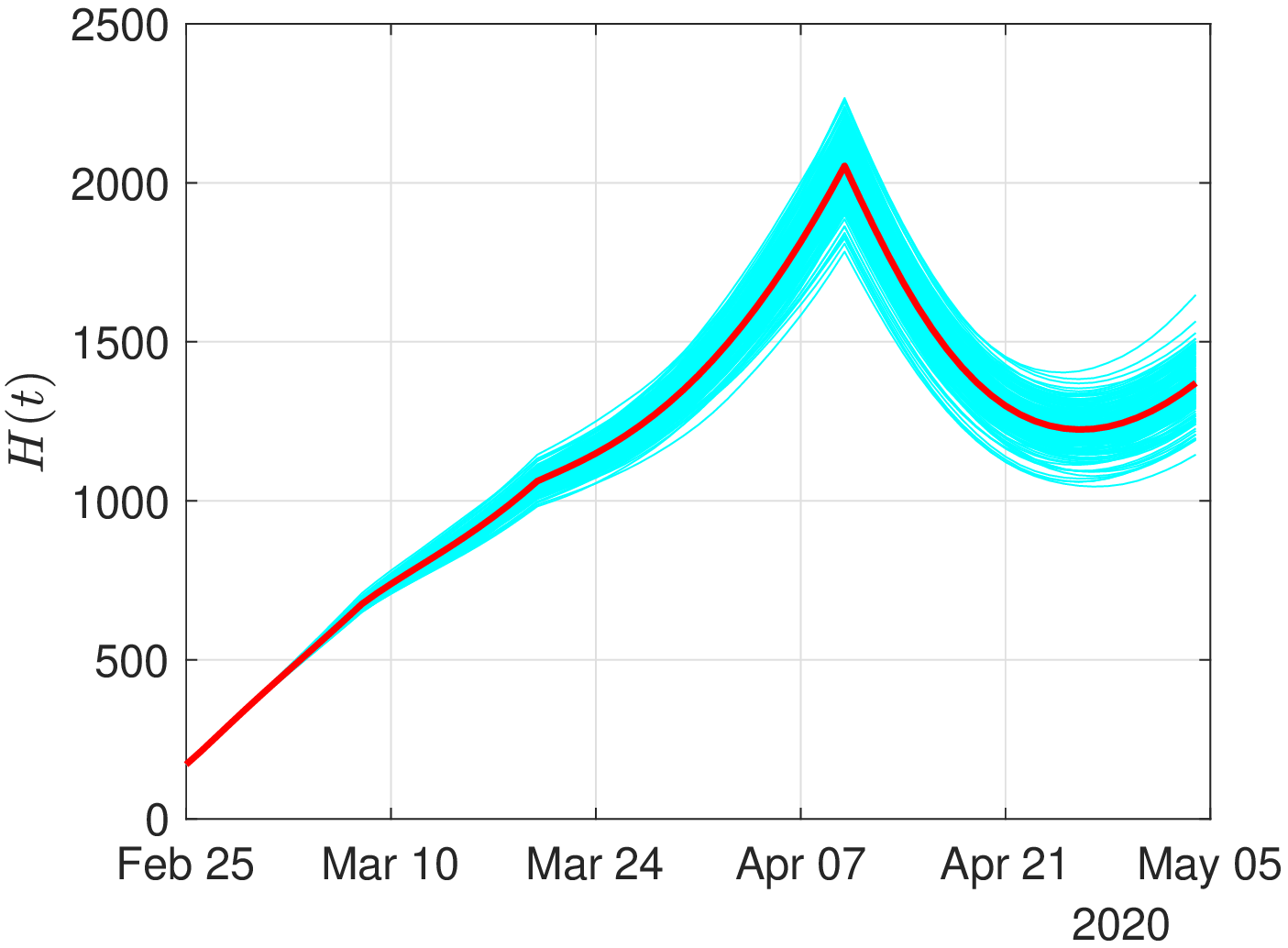}
    }
    \subfigure[]{
         \includegraphics[width=0.31\textheight]{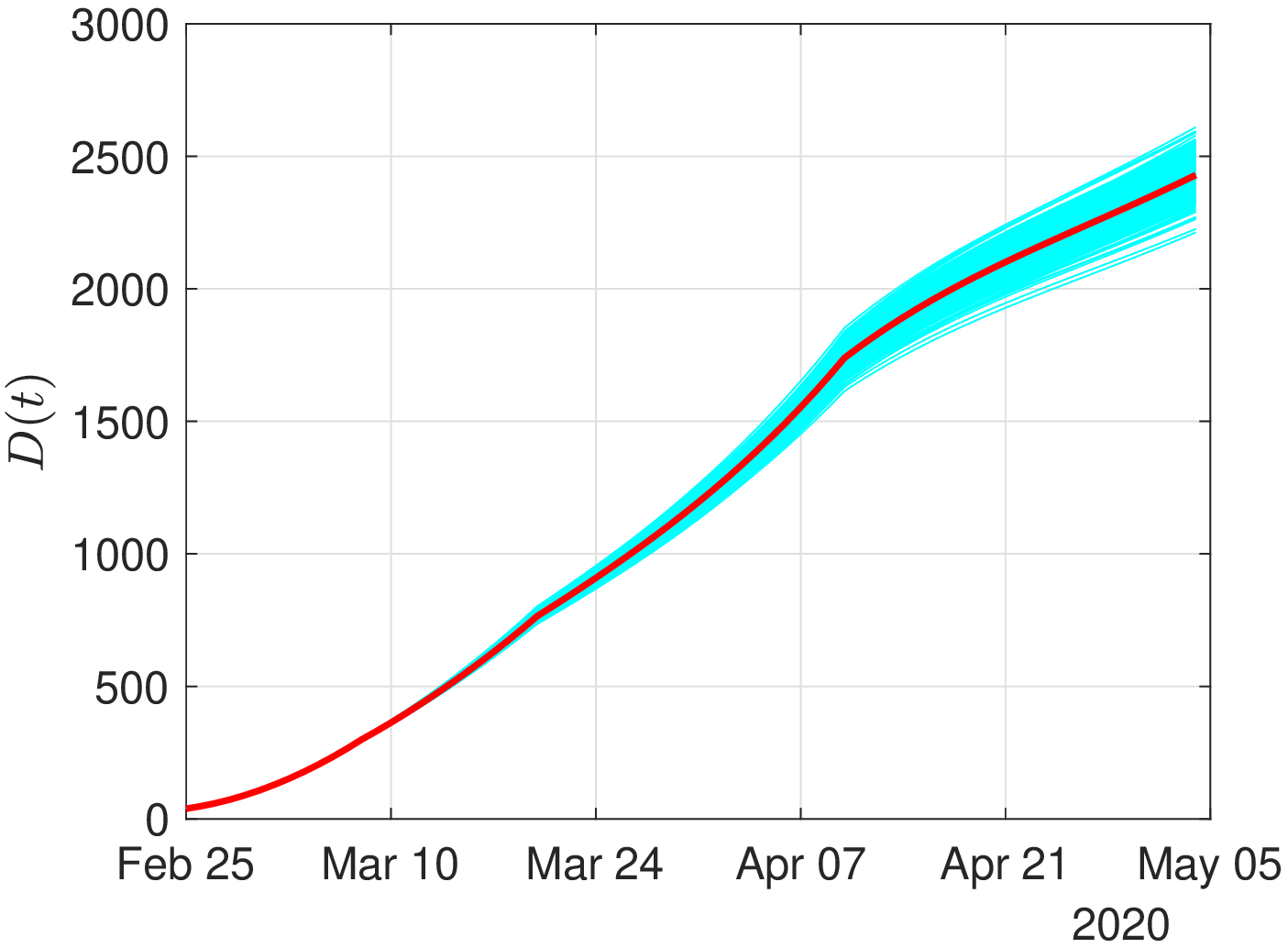}
    }
    \subfigure[]{
         \includegraphics[width=0.31\textheight]{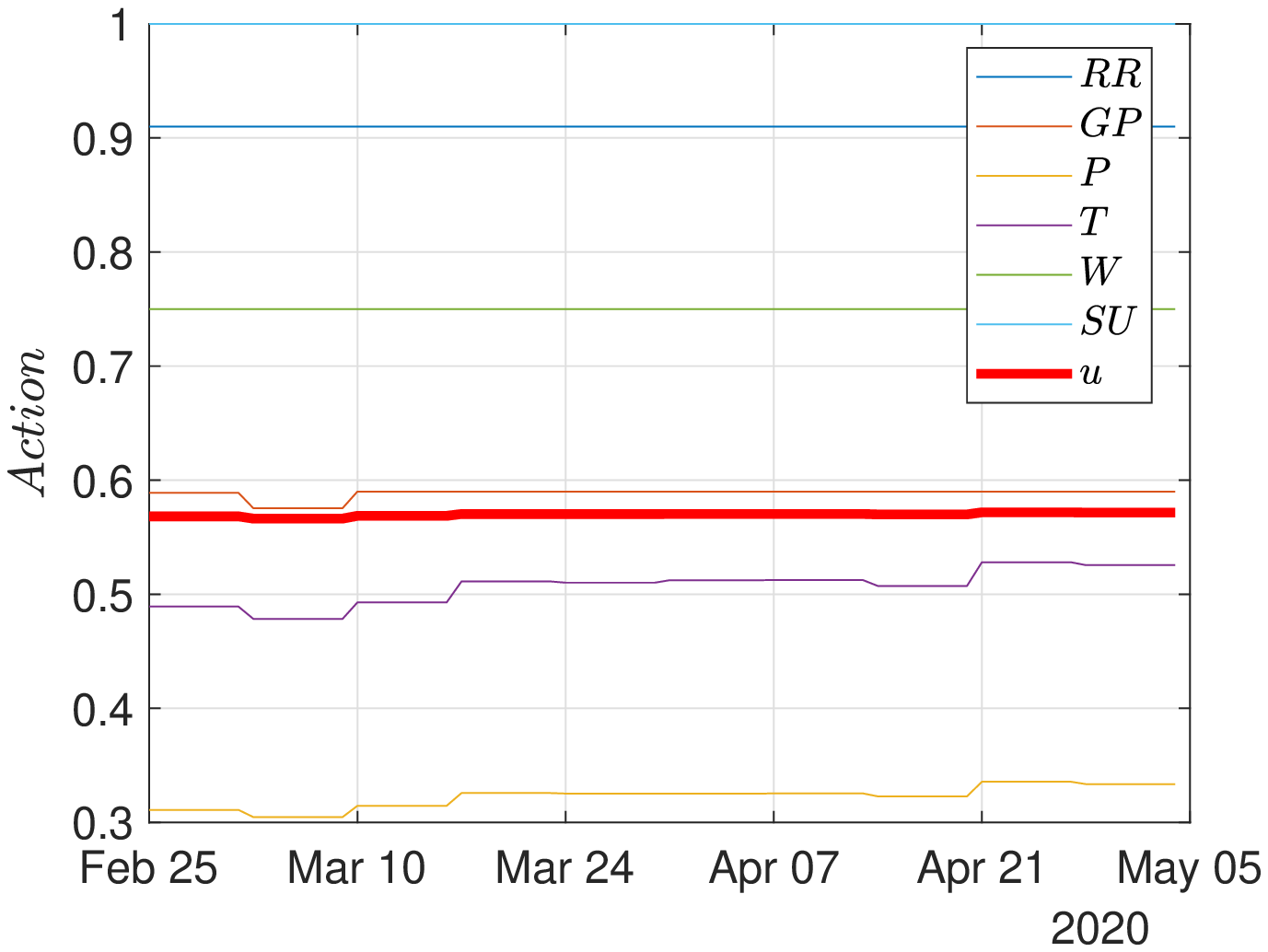}
    }
    \subfigure[]{
         \includegraphics[width=0.31\textheight]{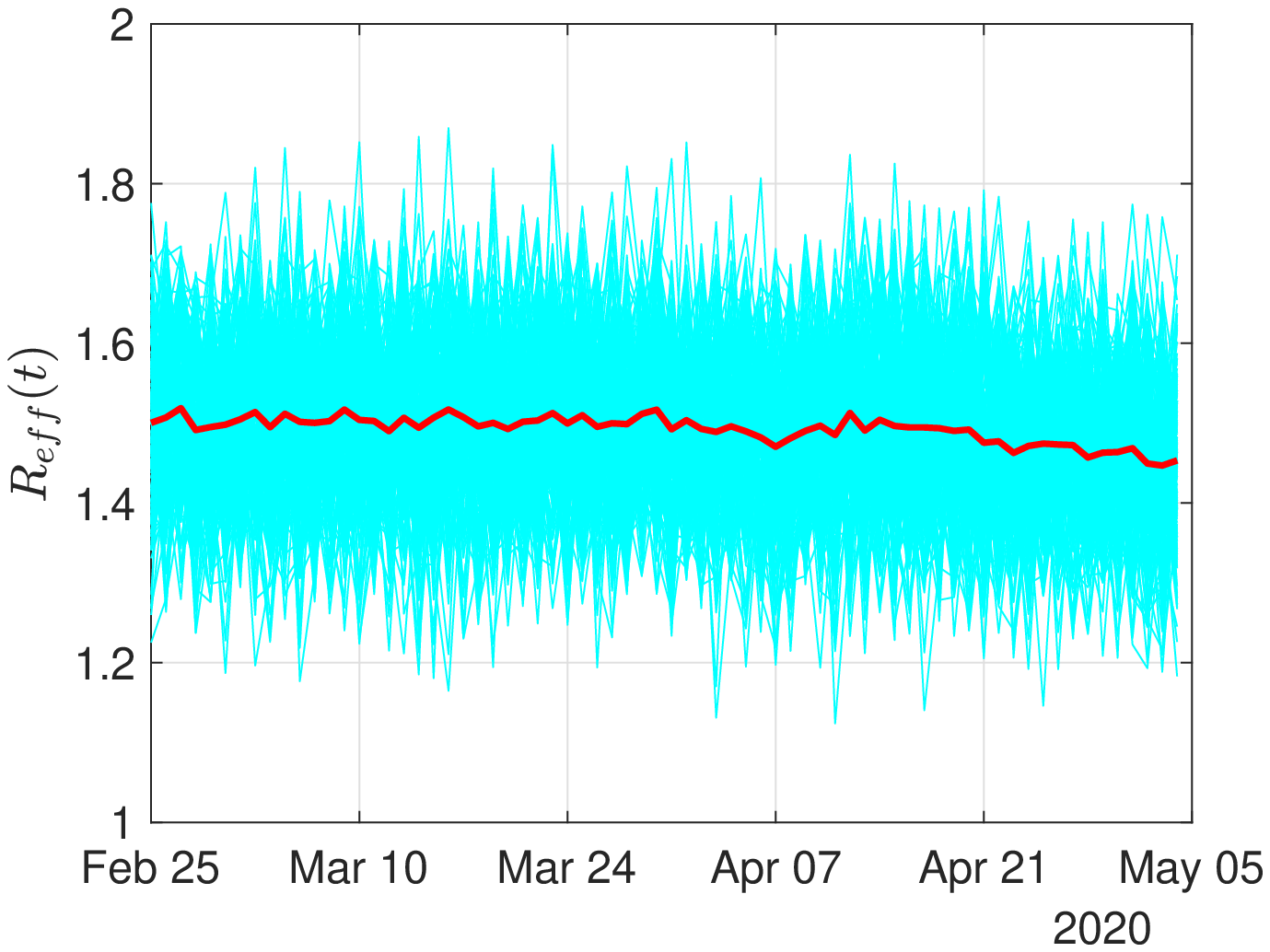}
    }
    \caption{Simulation results of epidemic progression under MPC for the period February 24-May 4 with $\omega_I=0$, $\omega_R=1$, $\Omega_A=0.6{\rm \bf I}_{6\times 6}$  and $\delta\alpha_c={\bf 1}_{6\times 1}$ in \eqref{OP1:0a}. (a) Infected cases. (b) Hospitalized. (c) Quarantined. (d) Dead. (e) Social distancing in the various activities resulting from the MPC. (f) Effective Reproduction number $R_e$. Red solid lines in (a), (b), (c), (d) \& (f) correspond to the trajectory expectation. The red solid line in (e) presents the cumulative mobility curtails. Cyan colors present individual trajectories of 200 adherence scenaria.}
    \label{f:MPC_S3}
\end{figure}

The effect of the choice of weights of $\Omega_A$ can be visualized for the case when $\omega_I=0$ \& $\omega_R=1$, i.e. they are relatively small compared to $\|\Omega_A\|_F=\approx 1.6$. We observe that when $\Omega_A=0.6{\rm \bf I}_{6\times 6}$
stringent curtails are imposed on important actives such as groceries and work (shown in Figure~\ref{f:MPC_S3}) compared to when the ``variable" weight is employed (shown in Figure~\ref{f:MPC_S2}), with the aggregate curtail being $u(t)\approx 0.57$. Considering the societal effects of activities the nominal weight for the optimization used in formulation of \eqref{OP1:J}-\eqref{OP1:3} was chosen as
\begin{equation}\label{e:V_A}
\Omega_A=\left[ 
    \begin{array}{cccccc}
0.2 & 0 & 0 & 0 & 0 & 0 \\
0 & 1.0 & 0 & 0 & 0 & 0 \\
0 & 0 & 0.5 & 0 & 0 & 0 \\
0 & 0 & 0 & 0.5 & 0 & 0 \\
0 & 0 & 0 & 0 & 1.0 & 0 \\
0 & 0 & 0 & 0 & 0 & 0.5 \\
\end{array}
\right]
\end{equation}

% \end{thebibliography}

\end{document}